\documentclass[12pt]{article}

\usepackage{doi, hyperref}
\usepackage{preprint-layout} %
\usepackage{preprint-notation} %

\usepackage{algorithmicx}
\usepackage{algpseudocode}
\usepackage{algorithm}
\algnewcommand\algorithmicforeach{\textbf{for each}}
\algdef{S}[FOR]{ForEach}[1]{\algorithmicforeach\ #1\ \algorithmicdo}

\renewcommand{\hatx}{\hat{x}}
\renewcommand{\haty}{\hat{y}}
\renewcommand{\hatr}{\hat{r}}

\newcommand{\hatO}{\widehat{\Omega}}

\newcommand{\hatnabla}{\widehat{\nabla}}

\newcommand{\hatV}{\widehat{V}}

\newcommand{\hatnu}{\widehat{\nu}}

\newcommand{\Reg}{R_{\delta}}
\newcommand{\Regi}{R_{i,\delta}}

\title{Robust Trimmed Multipatch IGA\\ with Singular Maps}
\author{Tobias Jonsson, Mats G. Larson, Karl Larsson}

\begin{document}

\maketitle

\begin{abstract}
We consider elliptic problems in multipatch isogeometric analysis (IGA) where the patch parameterizations may be singular. Specifically, we address cases where certain dimensions of the parametric geometry diminish as the singularity is approached -- for example, a curve collapsing into a point (in 2D), or a surface collapsing into a point or a curve (in 3D). To deal with this issue, we develop a robust weak formulation for the second-order Laplace equation that allows trimmed (cut) elements, enforces interface and Dirichlet conditions weakly, and does not depend on specially constructed approximation spaces. Our technique for dealing with the singular maps is based on the regularization of the Riemannian metric tensor, and we detail how to implement this robustly. We investigate the method's behavior when applied to a square-to-cusp parameterization that allows us to vary the singular behavior's aggressiveness in how quickly the measure tends to zero when the singularity is approached. We propose a scaling of the regularization parameter to obtain optimal order approximation. Our numerical experiments indicate that the method is robust also for quite aggressive singular parameterizations.
\end{abstract}

\section{Introduction}

Isogeometric analysis (IGA) combines the high fidelity geometry descriptions from CAD with powerful and versatile spline approximation spaces for generating accurate numerical solutions to PDE, see \cite{MR2152382,MR3618875}. 
The geometry in IGA is parameterized via a map
\begin{align}
F: \widehat{\Omega} \mapsto \Omega
\end{align}
from a reference domain $\widehat{\Omega} \subset \IR^{\hatd}$ to the physical domain $\Omega \subset\IR^{d}$, $d \geq \hatd$, where typically $d=2,3$.
Most often, $F$ is assumed to be regular, but in practice, singular parametrizations are not uncommon since they allow for flexible representations using a tensor product structure, for instance, mapping a square onto a triangle or a circle.

In this work, we mainly consider the type of singular parametrization where $F$ maps one part of the boundary of $\widehat{\Omega}$ onto a single point in the physical domain $\Omega$ or on its boundary $\partial\Omega$, but our proposed approach also applies to other types of singularities.
Our setting is a flexible multipatch framework, based on Nitsche's method for patch coupling, that also allows trimming of the reference domain, leading to so-called cut elements in the computational grid. In this framework, all terms in the weak form are formulated in the reference domain, where the effect of the map $F$ appears as the Riemannian metric tensor $G\in\IR^{ \widehat{d} \times \widehat{d}}$, which vary over the reference domain of each patch. When approaching a singular point, the matrix $G$ becomes singular, leading to divergent integrals in the weak form and numerical robustness issues.
To remedy this issue, we propose a general regularization procedure based on the eigenvalue decomposition of the metric tensor at each quadrature point.

\paragraph{Previous Works on Singular Maps.}
Even without any special treatment, IGA is often surprisingly robust concerning singular parametrizations, probably thanks to the use of Gauss quadrature, where it is unlikely that a quadrature point will be placed problematically close to a singular point \cite{MR3839189}. This is possibly why there is a comparably small number of IGA papers addressing this issue, but a series of contributions considering singular parameterizations are \cite{MR2851580, TAKACS2012361, MR3440759, 2015arXiv150708095T}. In \cite{MR2851580,TAKACS2012361}, a singular map of the same class addressed in the present work is included --- a polar parametrization of a circle that gives a singular point at the origin. It is shown that isogeometric spaces constructed using a standard approach will not belong to the suitable Sobolev space on $\Omega$. A correct subspace basis can be explicitly constructed by identifying a suitable subspace.
Approximation properties of such subspaces were analyzed in \cite{2015arXiv150708095T}.
The explicit construction of a correct subspace basis was generalized to correct subspaces of any smoothness in \cite{MR3440759} by composing the map as a singular map from a square onto a triangle, followed by a regular map onto the physical domain.
When using a standard isogeometric approximation space, i.e., not the correct subspace, numerical results indicate that the stability provided by the variational form lowers the continuity requirements of the approximation space. For instance, in \cite{MR3388842}, a polar parametrization of a sphere is considered, and it is demonstrated that enforcing only $C^0$ continuity at the poles is sufficient for optimal order convergence when solving a 4:th and a 6:th order PDE.

In contrast to most of these studies that aim to construct isogeometric bases with favorable properties, this work focuses on robustly computing method terms in the presence of singular parameterizations.

\paragraph{Trimmed Multipatch IGA.}

Large-scale CAD models of complex geometries rarely consist of single or matching NURBS patches. Instead, they are typically constructed using multiple trimmed patches connected along non-matching interfaces. In such a multipatch setting, it is natural to impose continuity over patch interfaces weakly, with the most common weak coupling approaches in IGA being penalty methods, mortar methods, and Nitsche-type methods, see \cite{MR3161053}.

Penalty methods \cite{MR351118,MR3353075} are popular due to their simplicity, but they have the drawback of aggressive $h$-scaling of the penalty term, which can result in poorly conditioned system matrices and requires careful tuning of the penalty parameter. Mortar methods \cite{MR3310287,MR3610098,MR3358117,MR4043895} introduce an auxiliary field over the interfaces, Lagrange multipliers, to ensure interface conditions are satisfied. This approach results in a saddle point problem, which is well-posed for multiplier field approximation spaces that fulfill the Ladyzhenskaya-Babuška-Brezzi (LBB) condition. However, mortar methods typically perform satisfactorily even without strictly meeting this condition. Nitsche-type methods \cite{MR341903,MR3144632,MR3510022,MR3630844} avoid introducing Lagrange multipliers and, unlike penalty methods, provide a consistent formulation that is well-conditioned even for higher-order polynomials. In this work, we utilize a Nitsche-type formulation for coupling the patches.

Since patches can be trimmed, the patch domains may cut arbitrarily through the computational grids, potentially causing ill-conditioning of the discretized problem. To ensure robustness regardless of the cut situation, stabilization is required. Stabilization techniques have been a research focus for many unfitted finite element methods, particularly cut finite element methods (CutFEM) \cite{MR3416285}, which typically employs so-called ghost penalty stabilization, introduced in \cite{MR2738930} and generalized to higher-order in \cite{MR3268662}. Other viable options for stabilizing trimmed Nitsche-type methods include discrete extension \cite{2022arXiv220506543B}, least-squares stabilization \cite{MR3920978,MR4394710}, basis function removal \cite{ElfLarLar18}, or cell merging \cite{MR3032318,MR3795012}. In this work, we use ghost penalty stabilization.

In special cases we might have precise control over how the multipatch geometry is generated, making it possible to construct matching approximation spaces with strongly enforced continuity. This results in a globally conforming approximation space, eliminating the need for additional weak coupling terms or stabilization due to trimming \cite{MR4489838}. However, such a construction is generally a non-trivial problem, especially if $C^1$-continuity is desired \cite{MR4254141,MR4673273}. A class of multipatch spline constructions relevant to the regularization procedure in this work is unstructured splines, as some are based on singular parameterizations, for instance the D-patch \cite{MR1462260}.

\paragraph{Challenges.}
The setting for this work is the framework developed in \cite{MR3682761,MR3709202,MR3959738,gap2022} for isogeometric methods on trimmed multipatch surfaces, using a Nitsche-type approach to couple the patches. While this framework is not inherently limited to surfaces, the focus on surfaces embedded in $\IR^3$ made a Riemann geometry approach natural for formulating the methods in the reference domain, as we do in the present work.
Two challenges of devising a method in such a framework that is robust with respect to singular maps are: 
\begin{itemize}
\item
The decoupling between the geometry description and the approximation space allows for trimmed reference domains and cut elements in the approximation space but, on the other hand, makes it awkward to employ techniques where much is assumed to be known regarding the geometry, such as explicitly constructing a correct subspace basis.

\item
The presence of boundary terms due to the weak enforcement of interface and boundary conditions makes it less likely that a naive approach, assuming that quadrature points will not be placed problematically close to a singularity, will work robustly.
\end{itemize}

\paragraph{Contributions.} The main contributions of this work are:
\begin{itemize}
	\item We propose a general procedure for dealing with singular parameterizations based on regularizing the Riemannian metric tensor. We apply this procedure to a multipatch IGA framework, allowing for trimmed isogeometric approximation spaces.
	\item We detail how to robustly compute the regularized quantity without running into numerical accuracy issues. Also when no regularization is desired, this appears to be a more numerically robust way of computing the terms involved.
	\item We sketch out a proof of how the consistency error induced by the regularization depends on the regularization parameter $\delta$ and present how $\delta$ should be scaled by the mesh size $h$ for an optimal order method. This is confirmed in numerical experiments.
\end{itemize}
In addition, we emphasize the following observations regarding the resulting method:
\begin{itemize} 
\item The proposed regularization procedure is efficient since it essentially requires no prior information about the singularity to be applied, which is particularly useful in the flexible multipatch IGA framework we work in.
\item The regularized multipatch formulation is guaranteed stable. It features the correct scaling of the Nitsche penalty parameter with respect to the Riemannian metric tensor and the mesh size of the joining patches.
\item The cost of the regularization is solving the very small eigenvalue problem $G a = \lambda a$ at each quadrature point, but, on the other hand, there is no longer a need to compute the inverse $G^{-1}$ at each quadrature point.
\item Numerical experiments indicate good performance and robustness, even when applied to very aggressive singularities, and the existence of functions in the approximation space outside the correct subspace does not seem to be an issue.
\end{itemize}

\paragraph{Outline.} The remainder of this paper is dispositioned as follows. In Section~\ref{sec:geometry}, we define our multipatch geometry and detail a model singular parameterization. In Section~\ref{sec:method}, we present a multipatch method and propose a regularization of the Riemannian metric tensor for dealing with the singular map. In Section~\ref{sec:numerics} we present numerical experiments and illustrating examples. Finally, in Section~\ref{sec:conclusions}, we give some summarizing conclusions.

\section{Multipatch Geometry}
\label{sec:geometry}

\subsection{Parametric Multipatch Geometry}

\begin{figure}
	\centering
	\includegraphics[width=0.70\linewidth]{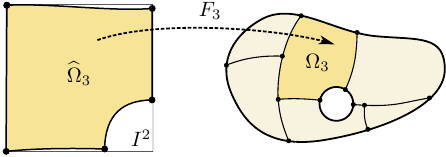}
	\caption{\emph{Parametric multipatch domain.} Multipatch domain consisting of six patches. Each patch $\Omega_i$ is parametrically described via a parametric map $F_i:\widehat{\Omega}_i \to \Omega_i$, from the patch reference domain $\widehat{\Omega}_i$ onto the patch $\Omega_i$. Note that patch reference domain may be trimmed and that, since the solution is coupled weakly over the patch interfaces, the patch decomposition allows for T-junctions.}
	\label{fig:multipatch}
\end{figure}

We first define our multipatch domain. Let $\Omega \subset \IR^d$ be a union of disjoint patches $\Omega_i$
\begin{equation}
\Omega = \bigcup_{i\in\mcI} \Omega_i
\end{equation}
where the boundary $\partial \Omega_i$ of each patch is piecewise smooth.
The interface between two patches $\Omega_i$ and $\Omega_j$, $i\neq j$, is expressed
\begin{align}
\partial\Omega_i \cap \partial\Omega_j
\end{align}
and the domain of all interfaces in $\Omega$ is
\begin{align}
\Gamma = \bigcup_{\substack{i,j\in\mcI \\ i\neq j}} \partial\Omega_i \cap \partial\Omega_j
\end{align}
The boundary of $\Omega$ is given by
\begin{align}
\partial\Omega = \Bigl( \bigcup_{i\in\mcI} \partial\Omega_i \Bigr) \setminus \Gamma \,,
\end{align}
i.e., the parts of the patch boundaries that are not interfaces.
Each patch $\Omega_i$ is described via a surjective parametric map $F_i:\widehat{\Omega}_i \to \Omega_i$ such that
\begin{align}
\Omega_i = F_i(\widehat{\Omega}_i)
\end{align}
where $\widehat{\Omega}_i \subset \IR^{\widehat{d}}$, $\widehat{d}\leq d$, is the reference domain, which typically is a subset of the unit cube in $\IR^{\widehat{d}}$.
This description is illustrated in Figure~\ref{fig:multipatch}.

The Jacobian of $F_i$ in $\hatO$ is
\begin{equation} \label{eq:mapderivatives}
DF_i = F_i\otimes\hatnabla
\end{equation}
where $\hatnabla$ is the gradient with respect to reference domain coordinates, and the Riemannian metric tensor $G_i$ in $\hatO_i$ is given by
\begin{align}
G_i = DF_i^T DF_i
\end{align}
When the patch association is evident from the context, we sometimes drop the subscript $i$ to simplify notation.

\subsection{Model for Singular Parametrization} \label{sec:model-singular-map}

\begin{figure}
\centering
\begin{subfigure}[t]{0.23\linewidth}\centering	
\includegraphics[width=0.95\linewidth]{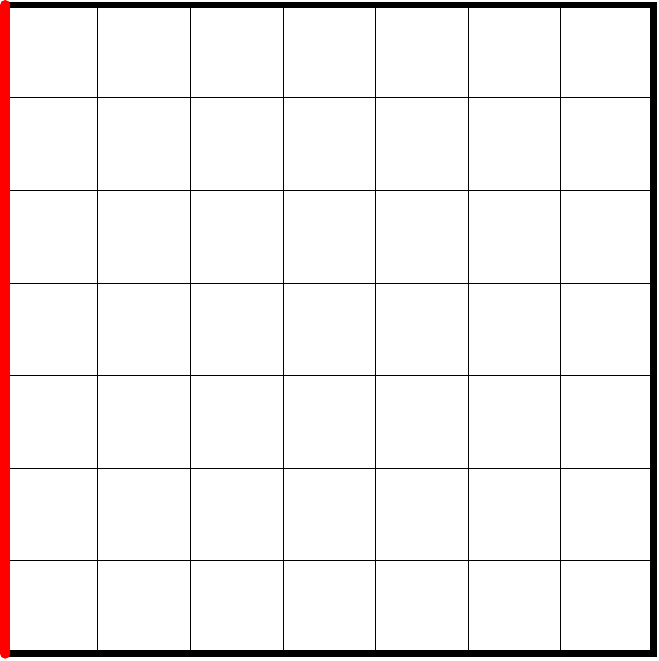}
\subcaption{$\widehat{\Omega}$}
\label{fig:map:1}
\end{subfigure}
\begin{subfigure}[t]{0.23\linewidth}\centering
\includegraphics[width=0.95\linewidth]{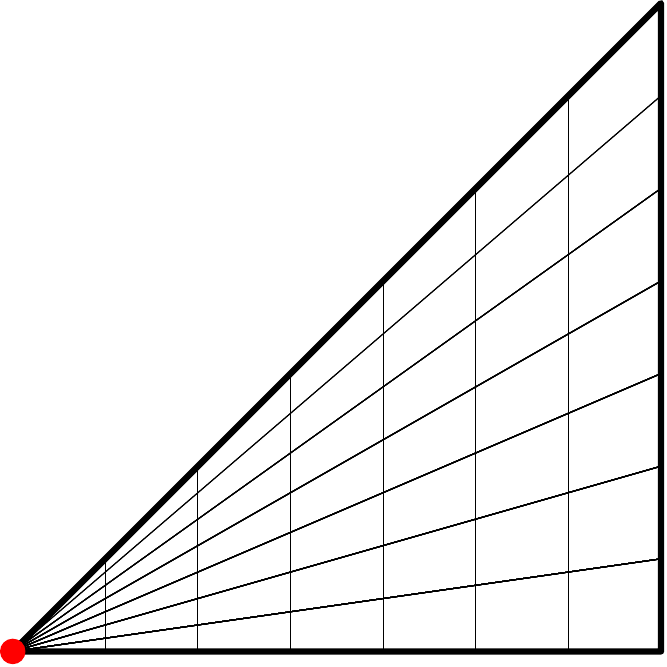}
\subcaption{$\Omega$, $\gamma=1$}
\label{fig:map:2}
\end{subfigure}
\begin{subfigure}[t]{0.23\linewidth}\centering
\includegraphics[width=0.95\linewidth]{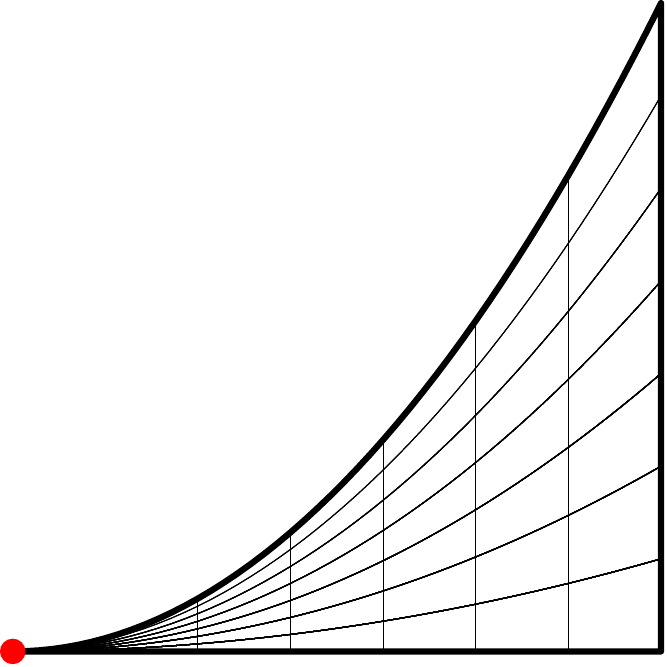}
\subcaption{$\Omega$, $\gamma=2$}
\label{fig:map:3}
\end{subfigure}
\begin{subfigure}[t]{0.23\linewidth}\centering	
\includegraphics[width=0.95\linewidth]{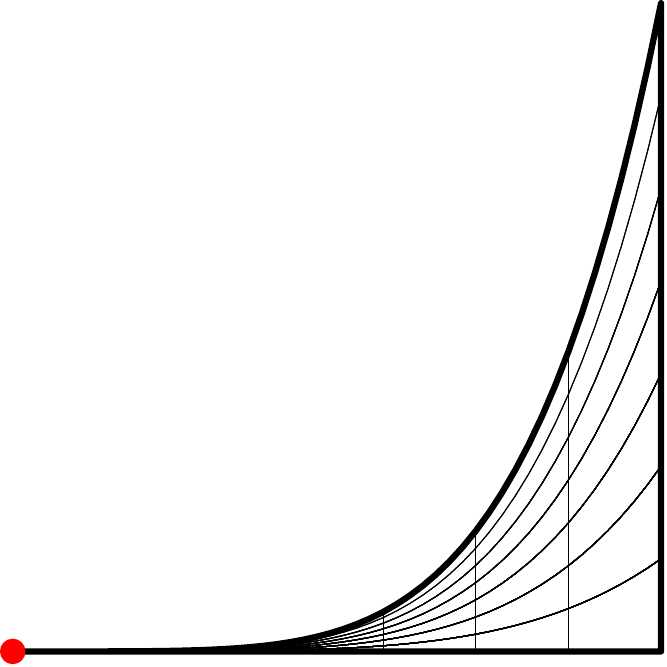}
\subcaption{$\Omega$, $\gamma=5$}
\label{fig:map:4}
\end{subfigure}
\caption{\emph{Square to cusp parameterization.} The leftmost figure shows the unit square reference domain equipped with a structured grid. The other figures show the effect of the map on this grid for various choices of $\gamma$.}
\label{fig:cuspmap}
\end{figure}

As a model for a singular map, we consider the following two-dimensional square-to-cusp parameterization.
Let $\hatO$ be the unit square and 
$\Omega = F(\hatO)$ with $F$ as the mapping
\begin{equation} \label{eq:square-to-cusp}
F
\left(
\left[
\begin{matrix}
\hatx
\\
\haty
\end{matrix}
\right]
\right)
=
\left[
\begin{matrix}
\hatx
\\
\hatx^\gamma \haty 
\end{matrix}
\right]
\end{equation}
where $\gamma\geq 1$ is a parameter that allows us to vary the aggressiveness of the cusp in the sense that the higher the value of $\gamma$, the faster the area measure $|G|^{1/2}$ vanishes when $\hatx\to 0$. See the illustration in Figure~\ref{fig:cuspmap} for examples with various choices of $\gamma$.
The Jacobian of $F$ in $\hatO$ is
\begin{equation}
DF =
\begin{bmatrix}
\frac{\partial F}{\partial \hatx} & \frac{\partial F}{\partial \haty}
\end{bmatrix}
=
\left[
\begin{matrix}
1 & 0 
\\
\gamma \hatx^{\gamma - 1} \haty & \hatx^\gamma
\end{matrix}
\right]
\end{equation}
and the Riemannian metric tensor $G$ in $\hatO$ is
\begin{equation}
G = DF^T DF = \left[
\begin{matrix}
1 & \gamma \hatx^{\gamma - 1}  \haty
\\
0 & \hatx^\gamma
\end{matrix}
\right]
\,
\left[
\begin{matrix}
1 & 0 
\\
\gamma \hatx^{\gamma - 1}  \haty & \hatx^\gamma
\end{matrix}
\right]
= 
\left[
\begin{matrix}
1 +\gamma^2 \hatx^{2(\gamma - 1)} \haty^2 & \gamma \hatx^{2\gamma-1} \haty
\\
 \gamma \hatx^{2\gamma-1} \haty & \hatx^{2\gamma} 
\end{matrix}
\right]
\end{equation}
with determinant
\begin{align}
|G|
=
(1 +\gamma^2 \hatx^{2(\gamma - 1)} \haty^2)  \hatx^{2\gamma}  
-
 \gamma^2 \hatx^{2(2\gamma-1)} \haty^2
 = \hatx^{2\gamma}
\end{align}
The eigenvalues to $G$ solves 
\begin{align}
0&= 
\left|
\begin{matrix}
1 +\gamma^2 \hatx^{2(\gamma - 1)} \haty^2 -\lambda &  \gamma \hatx^{2\gamma-1} \haty
\\
\gamma \hatx^{2\gamma-1} \haty & \hatx^{2\gamma}   -\lambda
\end{matrix}
\right|
\\&=
(1 +\gamma^2 \hatx^{2(\gamma - 1)} \haty^2-\lambda)(\hatx^{2\gamma}   -\lambda) - \gamma^2 \hatx^{2(2\gamma-1)} \haty^2
\\
&=
\lambda^2 
- ( 1+ \hatx^{2\gamma} + \gamma^2 \hatx^{2(\gamma-1)} \haty^2 ) \lambda 
+ \hatx^{2\gamma}
\end{align}
which gives
\begin{equation}
\lambda = \frac{1+ \hatx^{2\gamma} + \gamma^2 \hatx^{2(\gamma-1)} \haty^2 }{2} \pm \left(  \left( \frac{1+ \hatx^{2\gamma} + \gamma^2 \hatx^{2(\gamma-1)} \haty^2 }{2}\right)^2 - \hatx^{2\gamma}  \right)^{1/2}
\end{equation}
For small $\hatx$ we have the approximation 
\begin{align}
&\left(  \left( \frac{1+ \hatx^{2\gamma} + \gamma^2 \hatx^{2(\gamma-1)} \haty^2 }{2}\right)^2 - \hatx^{2\gamma}  \right)^{1/2}
\\
&\qquad\approx
\left( \left( \frac{1+ \hatx^{2\gamma} + \gamma^2 \hatx^{2(\gamma-1)} \haty^2 }{2}\right)^2   \right)^{1/2}
\\&\qquad\quad
-
\frac{1}{2} \left( \left( \frac{1+ \hatx^{2\gamma} + \gamma^2 \hatx^{2(\gamma-1)} \haty^2 }{2}\right)^2 \right)^{-1/2} \hatx^{2\gamma}
\\
&\qquad =
\frac{1+ \hatx^{2\gamma} + \gamma^2 \hatx^{2(\gamma-1)} \haty^2 }{2}
-
\frac{1}{2} \left( \frac{1+ \hatx^{2\gamma} + \gamma^2 \hatx^{2(\gamma-1)} \haty^2 }{2}\right)^{-1} x^{2\gamma}
\end{align}
and we find that the leading terms are
\begin{equation}
\lambda_1 \sim 1,\qquad \lambda_2 \sim \hatx^{2\gamma}
\end{equation}
where we note that one eigenvalue vanishes as the left side of the unit square $\hatO$ is approached ($\hatx\to 0$).
We can also solve for the eigenvectors analytically with the eigenvector associated to $\lambda_2$ being
\begin{align} \label{eq:eigenvector}
a_2 =
\begin{bmatrix}
\frac{\hatx}{\gamma \haty}
\left( \frac{\lambda_2}{\hatx^{2\gamma}} - 1 \right)
\\
1
\end{bmatrix}
\sim
\begin{bmatrix}
\hatx^{2\gamma-1} \haty
\\
1
\end{bmatrix}
\end{align}
where we derived the leading order terms similar to the corresponding derivation for the eigenvalues.

This is an example of a type of singular parametrization where part of the reference domain boundary $\partial\hatO$ is mapped onto a single point in the physical space. At the end of Section~\ref{sec:brief-analysis} we also briefly consider another type of singular parameterization, where the singularity instead stems from the columns of the map derivatives \eqref{eq:mapderivatives} becoming colinear when approaching a certain point.

\section{Multipatch Isogeometric Method}
\label{sec:method}

\subsection{Model Interface Problem}

As a model multipatch problem we consider: for all $i\in\mcI$ find $u_i:\Omega_i\to\IR$ such that
\begin{alignat}{2}
- \Delta u_i &= f_i &\qquad &\text{in $\Omega_i$}
\label{eq:bulkpde}
\\ 
u_i - u_j &= 0 &\qquad &\text{on $\partial\Omega_i \cap \partial\Omega_j$}
\label{eq:interface1}
\\
\nabla_{n_i} u_i + \nabla_{n_j} u_j &=0 &\qquad &\text{on $\partial\Omega_i \cap \partial\Omega_j$}
\label{eq:interface2}
\\ 
u_i &= g_i &\qquad &\text{on $\partial\Omega_i \cap \partial\Omega$}
\label{eq:boundarycond}
\end{alignat}
for all $j\in\mcI\setminus\{j\}$,
where $f_i:\Omega_i\to\IR$ and $g_i:\partial\Omega_i \cap \partial\Omega\to\IR$ are given data.
Here, \eqref{eq:interface1} and \eqref{eq:interface2} are interface conditions and \eqref{eq:boundarycond} is a Dirichlet boundary condition.
On a standard (non-multipatch) domain, this is equivalent to the Dirichlet problem
\begin{equation}
-\Delta u = f \quad \text{in $\Omega$},
\qquad 
u = g \quad \text{on $\partial \Omega$}
\end{equation}

\subsection{Multipatch Method}

The method will be entirely formulated in the reference domain $\hatO_i \subset \IR^{\widehat{d}}$ of each patch, where we also discretize the problem. We employ a Riemann geometry approach that makes the formulation directly applicable also in cases where the patch $\Omega_i \subset \IR^d$ has a positive codimension, that is $\widehat{d}<d$. 
Any function in the physical domain $\Omega_i$ we can express in the reference domain $\hatO_i$ via the pull-back
\begin{align}
\hatv_i = v_i \circ F_i
\end{align}
and from here on, we let hats on functions denote precisely this.

\paragraph{Isogeometric Approximation Space.}

We equip the reference domain $\hatO_i$ of each patch with a mesh $\mcK_{h,i}$ of mesh size $h$ and equip this mesh with a finite element space $\hatV_{h,i}$. The isogeometric approximation space $V_{h,i}$ on $\Omega_i$ is then defined by the push-forward 
\begin{equation}
V_{h,i} = \{v:\Omega_i \rightarrow \IR \ | \ v \circ F_i = \hatv, \ \hatv \in \hatV_{h,i}\}
\end{equation}
and the approximation spaces on the whole multipatch domain $\Omega$ respectively on all the reference patches are defined as the direct sums
\begin{align}
V_h = \bigoplus_{i\in\mcI} V_{h,i}
\qquad\text{and}\qquad
\hatV_h = \bigoplus_{i\in\mcI} \hatV_{h,i}
\end{align}
This is essentially a standard isogeometric approximation space but without any incorporation of essential boundary or interface conditions since these will instead be enforced weakly. 
For now, we assume that each mesh $\mcK_{h,i}$ conforms to the patch geometry of $\hatO_i$ in the sense that $\mcK_{h,i}$ is a partition of $\hatO_i$. In Section~\ref{sec:stabilization} we will loosen this requirement such that it is sufficient that $\mcK_{h,i}$ covers $\hatO_i$.

\paragraph{Method Derivation.}

The equation \eqref{eq:bulkpde} posed in the reference domain $\hatO_i$ is 
\begin{align}
-|G_i|^{-1/2} \hatnabla \cdot \bigl( |G_i|^{1/2} G_i^{-1} \hatnabla \hatu_i \bigr)  = \hatf_i
\end{align}
Multiplying by a test function $\hatv_i$, integrating over the patch, and using a Green's formula, we obtain
\begin{align}
\int_{\hatO_i}  \hatf_i \hatv_i |G_i|^{1/2}
&=
-\int_{\hatO_i}  \bigl( |G_i|^{-1/2} \hatnabla \cdot \bigl( |G_i|^{1/2} G_i^{-1} \hatnabla \hatu_i \bigr)\bigr) \hatv_i |G_i|^{1/2}
\\
&=
-\int_{\hatO_i}  \bigl(\hatnabla \cdot  \bigl( |G_i|^{1/2} G_i^{-1} \hatnabla \hatu_i \bigr)\bigr) \hatv_i
\\
&=
\int_{\hatO_i} \bigl(|G_i|^{1/2} G_i^{-1}  \hatnabla \hatu_i \bigr) \cdot \hatnabla \hatv_i
-
\int_{\partial \hatO_i} \hatnu_i \cdot \bigl( |G_i|^{1/2} G_i^{-1} \hatnabla \hatu_i \bigr) \hatv_i
\label{eq:simplegreen}
\end{align}
where $\hatnu_i$ is the outward pointing normal to $\partial\hatO_i$.
These two terms are identical to the terms we would expect if performing this calculation in the physical domain, i.e.,
\begin{align}
\int_{\hatO_i} \bigl(|G_i|^{1/2} G_i^{-1}  \hatnabla \hatu_i \bigr) \cdot \hatnabla \hatv_i
&= \int_{\Omega_i} \nabla u_i \cdot \nabla v_i
\\
\int_{\partial \hatO_i} \hatnu_i \cdot \bigl( |G_i|^{1/2} G_i^{-1} \hatnabla \hatu_i \bigr) \hatv_i
&=
\int_{\partial \Omega_i} (n_i \cdot \nabla u_i) v_i
\end{align}
We define the set of points on each reference patch boundary that maps to the interfaces respectively to the domain boundary
\begin{align}
\partial \hatO_{i,\Gamma}
:= \bigl\{ \hatx \in \partial\hatO_i \, | \, F_i(\hatx) \in \Gamma  \bigr\}
,\quad
\partial \hatO_{i,\partial\Omega}
&:= \bigl\{ \hatx \in \partial\hatO_i \, | \, F_i(\hatx) \in \partial\Omega  \bigr\}
\end{align}
Note that, due to singular maps collapsing to single points, the remaining part of the reference domain boundary
\begin{align}
\partial \hatO_{i,\emptyset} := \partial\hatO_i \setminus (\partial \hatO_{i,\Gamma} \cup \partial \hatO_{i,\partial\Omega})
\end{align}
might be non-empty. However, such parts will not give any contribution to the boundary integral in \eqref{eq:simplegreen} since $\hatnu_i \cdot \bigl( |G_i|^{1/2} G_i^{-1} \bigr) = 0$ on $\partial \hatO_{i,\emptyset}$.

To formulate the method we define the following average operator on $\partial \hatO_{i,\Gamma} \cup \partial \hatO_{i,\partial\Omega}$
\begin{align} \label{eq:average}
\langle \hatv \rangle(\hatx)
=
\left\{
\begin{alignedat}{2}
&\kappa_{ij} \hatv_i(\hatx) + \kappa_{ji} \hatv_j(\hatx_j) &\quad &\text{on $\partial \hatO_{i,\Gamma}$}
\\
&0 &\quad &\text{on $\partial \hatO_{i,\partial\Omega}$}
\end{alignedat}
\right.
\end{align}
where $j\in\mcI\setminus\{i\}$ and $\hatx_j \in \partial\hatO_j$ are implicitly defined such that $\hatx$ and $\hatx_j$ map to the same point $F_i(\hatx) = F_j(\hatx_j) \in \Gamma$, and $\kappa_{ij}$ are weights $0<\kappa_{ij}=1-\kappa_{ji} < 1$.
In terms of this average, we can collect both the interface condition \eqref{eq:interface1} and the boundary condition \eqref{eq:boundarycond} in the patchwise condition
\begin{align} \label{eq:bc-and-ic}
\hatu_i - \langle \hatu \rangle &=
\left\{
\begin{alignedat}{2}
&0 &\qquad &\text{on $\partial\hatO_{i,\Gamma}$}
\\
&\hatg_i &\qquad &\text{on $\partial \hatO_{i,\partial\Omega}$}
\end{alignedat}
\right.
\end{align}
Summing \eqref{eq:simplegreen} over all patches gives
\begin{align} 
\sum_{i\in\mcI}
\int_{\hatO_i}  \hatf_i \hatv_i |G_i|^{1/2}
&=
\sum_{i\in\mcI}
\int_{\hatO_i} \bigl(|G_i|^{1/2} G_i^{-1}  \hatnabla \hatu_i \bigr) \cdot \hatnabla \hatv_i
\\&\qquad\quad
-
\int_{\partial \hatO_{i,\Gamma} \cup \partial \hatO_{i,\partial\Omega}} \hatnu_i \cdot \bigl( |G_i|^{1/2} G_i^{-1} \hatnabla \hatu_i \bigr) \hatv_i
\\
&=
\sum_{i\in\mcI}
\int_{\hatO_i} \bigl(|G_i|^{1/2} G_i^{-1}  \hatnabla \hatu_i \bigr) \cdot \hatnabla \hatv_i
\\&\qquad\quad
-
\int_{\partial \hatO_{i,\Gamma} \cup \partial \hatO_{i,\partial\Omega}} \hatnu_i \cdot \bigl( |G_i|^{1/2} G_i^{-1} \hatnabla \hatu_i \bigr) (\hatv_i - \langle \hatv \rangle)
\end{align}
where we in the last term can subtract $\langle \hatv \rangle$ without affecting the value since $\langle \hatv \rangle = 0$ on points mapping to the boundary $\partial\Omega$ and on points mapping to an interface, we utilize \eqref{eq:interface2}, that the sum of the fluxes over the interface is zero.
Further, by the patchwise condition \eqref{eq:bc-and-ic}, we, without affecting consistency, may append the following terms inside the sum on the right-hand side
\begin{align}
\label{eq:sym-term}
&-
\int_{\partial \hatO_{i,\Gamma} \cup \partial \hatO_{i,\partial\Omega}} (\hatu_i - \langle \hatu \rangle) \hatnu_i \cdot \bigl( |G_i|^{1/2} G_i^{-1} \hatnabla \hatv_i \bigr) 
\\
\label{eq:pen-term}
&+
\int_{\partial \hatO_{i,\Gamma} \cup \partial \hatO_{i,\partial\Omega}}\frac{\beta}{h} (\hatu_i - \langle \hatu \rangle) \hatnu_i \cdot \bigl( |G_i|^{1/2} G^{-1}_{i} \bigr)\hatnu_i (\hatv_i - \langle \hatv \rangle) 
\end{align}
along with the following terms inside the sum on the left-hand side
\begin{align}
&-
\int_{\partial \hatO_{i,\partial\Omega}} \hatg_i \hatnu_i \cdot \bigl( |G_i|^{1/2} G_i^{-1} \hatnabla \hatv_i \bigr) 
\\
&+
\int_{\partial \hatO_{i,\partial\Omega}}\frac{\beta}{h} \hatg_i \hatnu_i \cdot \bigl( |G_i|^{1/2} G^{-1}_{i} \bigr)\hatnu_i (\hatv_i - \langle \hatv \rangle) 
\end{align}
The term \eqref{eq:sym-term} gives us a symmetric form, and the term \eqref{eq:pen-term} is a penalty term where $\beta > 0$ is the Nitsche parameter.
The above derivation motivates the following method.

\paragraph{Standard Isogeometric Method.} Find $\hatu_h \in \hatV_h$ such that
\begin{align}
\hata_h(\hatu_h,\hatv) = \hatl_h(\hatv) \qquad\text{for all $\hatv\in\hatV_h$}
\end{align}
with forms
\begin{align}
\hata_h(\hatv,\hatw)
&=
\sum_{i\in\mcI}
\int_{\hatO_i} \bigl(|G_i|^{1/2} G_i^{-1}  \hatnabla \hatv_i \bigr) \cdot \hatnabla \hatw_i
\\&\qquad\quad
-
\int_{\partial \hatO_{i,\Gamma} \cup \partial \hatO_{i,\partial\Omega}} \hatnu_i \cdot \bigl( |G_i|^{1/2} G_i^{-1} \hatnabla \hatv_i \bigr) (\hatw_i - \langle \hatw \rangle)
\\&\qquad\quad
-
\int_{\partial \hatO_{i,\Gamma} \cup \partial \hatO_{i,\partial\Omega}} (\hatv_i - \langle \hatv \rangle) \hatnu_i \cdot \bigl( |G_i|^{1/2} G_i^{-1} \hatnabla \hatw_i \bigr) 
\\&\qquad\quad
+
\int_{\partial \hatO_{i,\Gamma} \cup \partial \hatO_{i,\partial\Omega}}\frac{\beta}{h} (\hatv_i - \langle \hatv \rangle) \hatnu_i \cdot \bigl( |G_i|^{1/2} G^{-1}_{i} \bigr)\hatnu_i (\hatw_i - \langle \hatw \rangle)
\end{align}
and
\begin{align}
\hatl_h(\hatw)
&=
\sum_{i\in\mcI}
\int_{\hatO_i}  \hatf_i \hatw_i |G_i|^{1/2}
\\&\qquad\quad-
\int_{\partial \hatO_{i,\partial\Omega}} \hatg_i \hatnu_i \cdot \bigl( |G_i|^{1/2} G_i^{-1} \hatnabla \hatw_i \bigr) 
\\&\qquad\quad+
\int_{\partial \hatO_{i,\partial\Omega}}\frac{\beta}{h} \hatg_i \hatnu_i \cdot \bigl( |G_i|^{1/2} G^{-1}_{i} \bigr)\hatnu_i (\hatw_i - \langle \hatw \rangle) 
\end{align}
A natural energy norm for this formulation is
\begin{align}
\tn \hatv \tn_h^2
&=
\sum_{i\in\mcI}
\bigl\|
|G_i|^{1/4} G_{i}^{-1/2} \hatnabla \hatv_i
\bigr\|_{L^2(\hatO_i)}^2
+
\bigl\|
h^{1/2} |G_i|^{1/4} G_{i}^{-1/2} \hatnabla \hatv_i
\bigr\|_{L^2({\partial \hatO_{i,\Gamma} \cup \partial \hatO_{i,\partial\Omega}})}^2
\\&\qquad\quad\nonumber +
\bigl\|
h^{-1/2} |G_i|^{1/4} G_{i}^{-1/2} \hatnu_i  (\hatv_i - \langle \hatv \rangle)
\bigr\|_{L^2({\partial \hatO_{i,\Gamma} \cup \partial \hatO_{i,\partial\Omega}})}^2
\end{align}
In a standard situation, this formulation gives a continuous and coercive bilinear form $\widehat{a}_h$, given that $\beta>0$ is sufficiently large.
However, since all terms in $\widehat{a}_h$ and $\tn \cdot \tn_h$ contain the expression $|G_i|^{1/2} G_i^{-1}$, depending on the map $F_i$, the integrals might be singular.

\begin{rem}[Nitsche Parameter Scaling]
By performing the entire derivation of the multipatch method in the reference domain, and formulating the Nitsche coupling terms patchwise via the average operator $\langle \cdot \rangle$, we directly obtain a Nitsche parameter that scales correctly with respect to the Riemannian metric tensor. This is of particular importance when considering maps that may be singular, since the metric tensor then may vary enormously, both along the patch boundaries and over the interfaces. By keeping the interface terms as patchwise expressions, rather than the common approach of collecting those terms interfacewise, we also directly obtain the correct weighting of the mesh sizes of the patches joined by an interface.
\end{rem}

\begin{rem}[Use of Standard Isogeometric Spaces]
	The fact that the expression $|G_i|^{1/2} G_i^{-1}$ is very large in regions close to singularities is likely the reason why a standard isogeometric approximation space, which can contain functions that don't belong to the correct Sobolev space, is feasible in a singular situation since this will force the solution to be essentially constant in certain directions within those regions.
\end{rem}

\subsection{Regularized Method}

For notational simplicity, we write this section for the case $\widehat{d}=2$ and assume one eigenvalue of the metric tensor $G$ tends to zero when approaching a singular point while the other remains bounded. Expanding $G$ in terms of eigenvectors $a_i$ and eigenvalues $\lambda_i$ we 
have 
\begin{align}
G = \sum_{k=1}^2 \lambda_k a_k \otimes a_k, 
\qquad 
G^{-1} = \sum_{k=1}^2 \lambda_k^{-1} a_k \otimes a_k
\end{align}
In order to regularize the problem we replace $\lambda_k$ by 
\begin{equation}
\lambda_{k,\delta}  = \lambda_k^{1/2}  \max (\lambda_k^{1/2},\delta^{1/2})
\end{equation}
where $\delta$ is a positive parameter
and we define the regularized metric tensor
\begin{equation}
G_\delta = \sum_{k=1}^2 \lambda_{k,\delta} a_k \otimes a_k
\end{equation}
We then have 
\begin{align}
|G|^{1/2} G_\delta^{-1}
&= 
 \sum_{k=1}^2 ( \lambda_1 \lambda_2 )^{1/2} \lambda_{k,\delta} ^{-1} a_k \otimes a_k 
\\ \label{eq:G-scaling-identity}
& = 
\frac{\lambda_2^{1/2}} {\max (\lambda_1^{1/2},\delta^{1/2})}    a_1 \otimes a_1 
+
\frac{\lambda_1^{1/2}} {\max (\lambda_2^{1/2},\delta^{1/2})}   a_2 \otimes a_2 
\\
& \leq
\frac{\lambda_2^{1/2}}{\delta^{1/2}} a_1 \otimes a_1 
 +
  \frac{\lambda_1^{1/2}}{\delta^{1/2} } a_2 \otimes a_2 
\\ \label{eq:Gtilde-bounded}
&< \infty
\end{align}
for $\delta>0$. By replacing $|G|^{1/2} G^{-1}$ in the standard method by \eqref{eq:G-scaling-identity} we get the following regularized version of the method.
How to robustly compute the quantity corresponding to \eqref{eq:G-scaling-identity} when $\widehat{d}\geq 2$ is provided in Algorithm~\ref{alg:identity}.

\paragraph{Regularized Isogeometric Method.} Find $\hatu_h \in \hatV_h$ such that
\begin{align}
\hata_{h,\delta}(\hatu_h,\hatv) = \hatl_{h,\delta}(\hatv) \qquad\text{for all $\hatv\in\hatV_h$}
\end{align}
with forms
\begin{align}
\hata_{h,\delta}(\hatv,\hatw)
&=
\sum_{i\in\mcI}
\int_{\hatO_i} \bigl(\Regi  \hatnabla \hatv_i \bigr) \cdot \hatnabla \hatw_i
\\&\qquad\quad
-
\int_{\partial \hatO_{i,\Gamma} \cup \partial \hatO_{i,\partial\Omega}} \hatnu_i \cdot \bigl( \Regi \hatnabla \hatv_i \bigr) (\hatw_i - \langle \hatw \rangle)
\\&\qquad\quad
-
\int_{\partial \hatO_{i,\Gamma} \cup \partial \hatO_{i,\partial\Omega}} (\hatv_i - \langle \hatv \rangle) \hatnu_i \cdot \bigl( \Regi \hatnabla \hatw_i \bigr) 
\\&\qquad\quad
+
\int_{\partial \hatO_{i,\Gamma} \cup \partial \hatO_{i,\partial\Omega}}\frac{\beta}{h} (\hatv_i - \langle \hatv \rangle) \hatnu_i \cdot \Regi \hatnu_i (\hatw_i - \langle \hatw \rangle)
\end{align}
and
\begin{align}
\hatl_{h,\delta}(\hatw)
&=
\sum_{i\in\mcI}
\int_{\hatO_i}  \hatf_i \hatw_i |G_i|^{1/2}
\\&\qquad\quad-
\int_{\partial \hatO_{i,\partial\Omega}} \hatg_i \hatnu_i \cdot \bigl( \Regi \hatnabla \hatw_i \bigr) 
\\&\qquad\quad+
\int_{\partial \hatO_{i,\partial\Omega}}\frac{\beta}{h} \hatg_i \hatnu_i \cdot  \Regi \hatnu_i (\hatw_i - \langle \hatw \rangle) 
\end{align}
where
\begin{align}
\Regi = |G_i|^{1/2} G_{i,\delta}^{-1}
\end{align}
is computed according to the identity \eqref{eq:G-scaling-identity} or, more generally, using Algorithm~\ref{alg:identity}.
It is important to use this identity for the regularization to work.

\begin{algorithm}
\textbf{Input:} A quadrature point $\hatx \subset \overline{\widehat{\Omega}}$\\
\textbf{Output:} $\Reg = |G|^{1/2} G^{-1}_{\delta}$
\begin{algorithmic}[1]
\State EVP: Find all pairs $\{(\lambda_k,a_k)\}_{k=1}^{\widehat{d}}$ such that $G(\widehat{x})a_k = \lambda_k a_k$
\State Initialize: $\Reg = 0$
\ForEach {$(\lambda_k,a_k)$}
\State  $\Reg \gets \Reg +
\bigl(\prod_{j=1, j\neq k}^{\widehat{d}} \lambda_j^{1/2} \bigr)
\bigl( (a_k\cdot a_k) \max\bigl(\delta^{1/2},\lambda_k^{1/2} \bigr) \bigr)^{-1} a_k \otimes a_k$ 
\EndFor
\end{algorithmic}
\caption{Robust computation of $\Reg = |G|^{1/2} G^{-1}_{\delta}$.}
\label{alg:identity}
\end{algorithm}

A natural energy norm for the regularized formulation is
\begin{align} \label{eq:energy-norm-reg}
\tn \hatv \tn_{h,\delta}^2
&=
\sum_{i\in\mcI}
\bigl\|
	\Reg^{1/2} \hatnabla \hatv_i
\bigr\|_{L^2(\hatO_i)}^2
+
\bigl\|
h^{1/2} \Reg^{1/2} \hatnabla \hatv_i
\bigr\|_{L^2({\partial \hatO_{i,\Gamma} \cup \partial \hatO_{i,\partial\Omega}})}^2
\\&\qquad\quad +
\bigl\|
h^{-1/2} \Reg^{1/2} \hatnu_i  (\hatv_i - \langle \hatv \rangle)
\bigr\|_{L^2({\partial \hatO_{i,\Gamma} \cup \partial \hatO_{i,\partial\Omega}})}^2
\end{align}
For $\delta>0$, we by the calculation \eqref{eq:G-scaling-identity}--\eqref{eq:Gtilde-bounded} have that $\Reg^{1/2}$ is bounded and, hence, the integrals in $\hata_{h,\delta}$ and $\tn \cdot \tn_{h,\delta}^2$ are not singular. Thereby, we have coercivity and continuity of $\hata_{h,\delta}$ in energy norm using standard arguments for Nitsche's method.

\begin{rem}[Applicability to Other PDEs] \label{rem:other-pde}
The numerical robustness of this regularization procedure stems from the key observation that the inverse metric tensor always appears in the form $|G|^{1/2}G^{-1}$, which allows us to handle cancellations between the measure and the inverse metric analytically, rather than numerically. Since this observation is problem specific, the corresponding term for other PDE must be identified analogously. For higher-order PDE, corresponding terms may also contain derivatives of the metric.
\end{rem}

\subsection{Trimmed Method}
\label{sec:stabilization}

We now loosen the requirement that the computational mesh $\mcK_{h,i}$ should be a partition of $\hatO_i$, and instead let it be sufficient that $\mcK_{h,i}$ covers $\hatO_i$, meaning that the patch boundary $\partial\hatO_i$ may arbitrarily cut through the computational mesh, creating elements that may have a very small intersection with the domain.
This is very convenient in IGA since it allows us to use approximation spaces with tensor product structure also when the geometric description includes trimming.
To maintain a coercive bilinear form in this setting, we require some extra stability close to the boundary of each patch. We provide this with the following so-called ghost penalty stabilization form
\begin{align}
\widehat{s}_{h,i}(\widehat{v}_i,\widehat{w}_i) = \eta \sum_{F \in \mathcal{F}_{h,i}}
\sum_{\ell=1}^p h^{2\ell-1} \int_F [\nabla^\ell_n \widehat{v}_i] [\nabla^\ell_n \widehat{w}_i]
\end{align}
where $\eta>0$ is a parameter, typically of a magnitude inversely proportional to the size of the Nitsche parameter $\beta$. Here, $\mathcal{F}_{h,i}$ is the set of all faces in the mesh $\mathcal{K}_{h,i}$ in the vicinity of the patch boundary $\partial\hatO_i$, and $[\nabla^\ell_n \widehat{v}_i]|_F$ is the jump in the $\ell$:th derivative of $\hatv_i$ in the normal direction of a face $F$.
By adding such stabilization, we arrive at the following method.

\paragraph{Regularized Trimmed Isogeometric Method.} Find $\hatu_h \in \hatV_h$ such that
\begin{align} \label{eq:the-method}
\hata_{h,\delta}(\hatu_h,\hatv) + \hats_{h}(\hatu_h,\hatv) = \hatl_{h,\delta}(\hatv) \qquad\text{for all $\hatv\in\hatV_h$}
\end{align}
where
\begin{align} \label{eq:ghost-penalty}
\widehat{s}_{h}(\widehat{v},\widehat{w})
=
\sum_{i\in\mcI}
\widehat{s}_{h,i}(\widehat{v}_i,\widehat{w}_i)
\end{align}
In this trimmed (cut) situation, we append the term $\widehat{s}_{h}(\widehat{v},\widehat{v})$ to the square of the energy norm \eqref{eq:energy-norm-reg} and can prove that we enjoy the same stability and accuracy properties as a fitted isogeometric method, see \cite{MR2738930,MR3682761}.

\subsection{Consistency Estimate}
\label{sec:brief-analysis}

As a motivation for how to choose the regularization parameter $\delta$ for an optimal order method, we here provide a sketch for an estimate of the consistency error induced by the regularization in the bulk term of the bilinear forms in a single patch $\Omega$.

\paragraph{General Arguments.}
In a Strang-type argument, the term associated with the consistency error in the bilinear form reads
\begin{align} \label{eq:strang-type}
\sup_{\hatv \in \hatV_h}
\frac{| \hata_h(\hatu,\hatv) - \hata_{h,\delta}(\hatu,\hatv)|}{\tn \hatv \tn_{h,\delta}}
\end{align}
where the square of the energy norm \eqref{eq:energy-norm-reg} includes the regularized bulk expression
\begin{align}
\int_{\hatO} (\Reg\hatnabla \hatv) \cdot \hatnabla \hatv
=
\int_{\hatO} (|G|^{1/2} G_{\delta }^{-1}\hatnabla \hatv) \cdot \hatnabla \hatv
=
\left\|
G_{\delta}^{-1/2} \hatnabla \hatv |G|^{1/4}
\right\|_{L^2(\hatO)}^2
\end{align}
Defining the region of regularization
\begin{align}
\hat{\omega}_\delta = \bigl\{ \hatx \in \hatO \, | \, \min_{k=1,\dots,\hatd} \lambda_k < \delta \bigr\}
\end{align}
and assuming that only the eigenvalue $\lambda_j$ needs to be regularized we have that
\begin{align}
G^{-1} - G_{\delta}^{-1}
&=
\sum_{k=1}^{\hatd}
\frac{1}{\lambda_k^{1/2}}
\left(\frac{1}{\lambda_k^{1/2}}-\frac{1}{\lambda_{k,\delta}^{1/2}}\right)a_k\otimes a_k
\\& \label{eq:Gdiff-simple}
=
\left\{
\begin{alignedat}{2}
&\frac{1}{\lambda_j^{1/2}}\left(\frac{1}{\lambda_j^{1/2}}-\frac{1}{\delta^{1/2}}\right)a_j\otimes a_j &\qquad &\text{in $\hat{\omega}_\delta$} \\
&0 &\qquad &\text{in $\hatO \setminus\hat{\omega}_\delta$}
\end{alignedat}
\right.
\end{align}
Considering the bulk part of the numerator in \eqref{eq:strang-type} and using \eqref{eq:Gdiff-simple} we obtain
\begin{align}
&\int_{\hatO} \bigl(|G|^{1/2} G^{-1}  \hatnabla \hatu \bigr) \cdot \hatnabla \hatv
-
\int_{\hatO} \bigl(\Reg  \hatnabla \hatu \bigr) \cdot \hatnabla \hatv
\\&\qquad
=
\int_{\hat{\omega}_\delta} \bigl(|G|^{1/2} ( G^{-1} - G_{\delta}^{-1} )   \hatnabla \hatu \bigr) \cdot \hatnabla \hatv
\\&\qquad
=
\int_{\hat{\omega}_\delta} \bigl(
G_{\delta}^{1/2} ( G^{-1} - G_{\delta}^{-1} ) G^{1/2} \lambda_j^{1/2} \lambda_j^{-1/2}G^{-1/2}\hatnabla \hatu  |G|^{1/4}\bigr) \cdot G_{\delta}^{-1/2} \hatnabla \hatv |G|^{1/4}
\\&\qquad
=
\int_{\hat{\omega}_\delta} \biggl(
\delta^{1/4} \biggl(1-\frac{\lambda_j^{1/2}}{\delta^{1/2}}\biggr)  a_j
\lambda_j^{-3/4} \partial_{a_j} \hatu  |G|^{1/4}\biggr) \cdot G_{\delta}^{-1/2} \hatnabla \hatv |G|^{1/4}
\\&\qquad
\leq \delta^{1/4}
\underbrace{
\Biggl\|
1-\frac{\lambda_j^{1/2}}{\delta^{1/2}}
\Biggr\|_{L^\infty(\hat{\omega}_\delta)}
}_{\sim 1}
\underbrace{
\bigl\| \lambda_j^{-1/2} \partial_{a_j} \hatu \bigr\|_{L^\infty(\hat{\omega}_\delta)}
}_{\leq \|\nabla u \|_{L^\infty(\omega_\delta)}}
\\&\qquad\qquad\quad \cdot 
\bigl\|
\lambda_j^{-1/4} |G|^{1/4}
\bigr\|_{L^2(\hat{\omega}_\delta)}
\underbrace{
\bigl\|
G_{\delta}^{-1/2} \hatnabla \hatv |G|^{1/4}
\bigr\|_{L^2(\hat{\omega}_\delta)}
}_{\leq \tn \hatv \tn_{h,\delta}}
\\&\qquad
\lesssim
\delta^{1/4}
\bigl\| \lambda_j^{-1/4} |G|^{1/4} \bigr\|_{L^2(\hat{\omega}_\delta)}
\|\nabla u \|_{L^\infty(\omega_\delta)}
\tn \hatv \tn_{h,\delta}
\label{eq:basic-est}
\end{align}
where the energy norm bound on the $\hatv$-term allows us to cancel the denominator in \eqref{eq:strang-type}, and we can use the max-norm bound on $\nabla u$ since $\sum_{k=1}^{\hatd} |\lambda_k^{-1/2} \partial_{a_k} \hatu |^2 = \| \nabla u \|_{\IR^d}^2$.
How the remaining factor $\bigl\| \lambda_j^{-1/4} |G|^{1/4} \bigr\|_{L^2(\hat{\omega}_\delta)}$ scales with $\delta$ depends on the specific map we are considering.

\begin{figure}
\centering
\begin{subfigure}[t]{0.23\linewidth}\centering	
\includegraphics[width=0.95\linewidth]{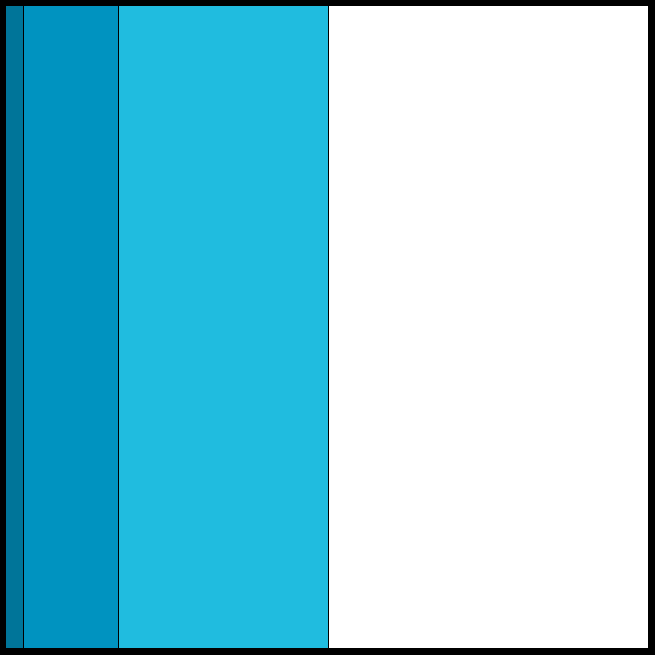}
\subcaption{$\hat{\omega}_\delta$, $\gamma=1,2,5$}
\label{fig:reg:1}
\end{subfigure}
\begin{subfigure}[t]{0.23\linewidth}\centering
\includegraphics[width=0.95\linewidth]{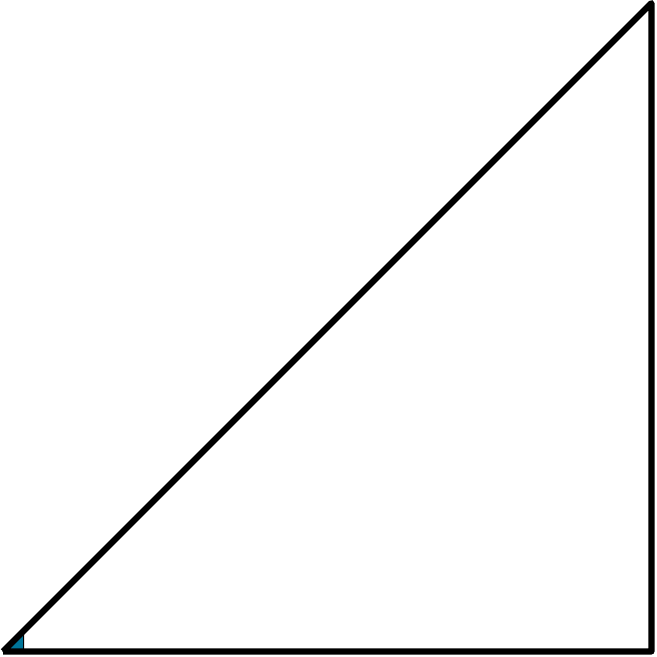}
\subcaption{$\omega_\delta$, $\gamma=1$}
\label{fig:reg:2}
\end{subfigure}
\begin{subfigure}[t]{0.23\linewidth}\centering
\includegraphics[width=0.95\linewidth]{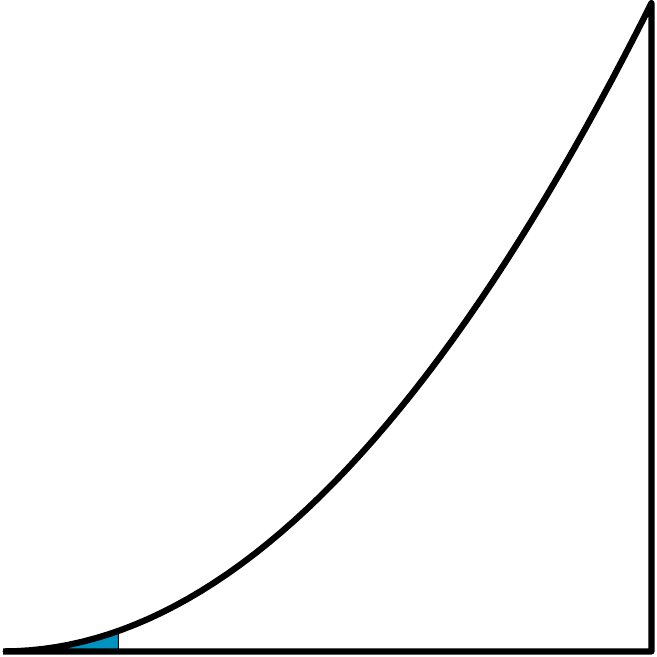}
\subcaption{$\omega_\delta$, $\gamma=2$}
\label{fig:reg:3}
\end{subfigure}
\begin{subfigure}[t]{0.23\linewidth}\centering	
\includegraphics[width=0.95\linewidth]{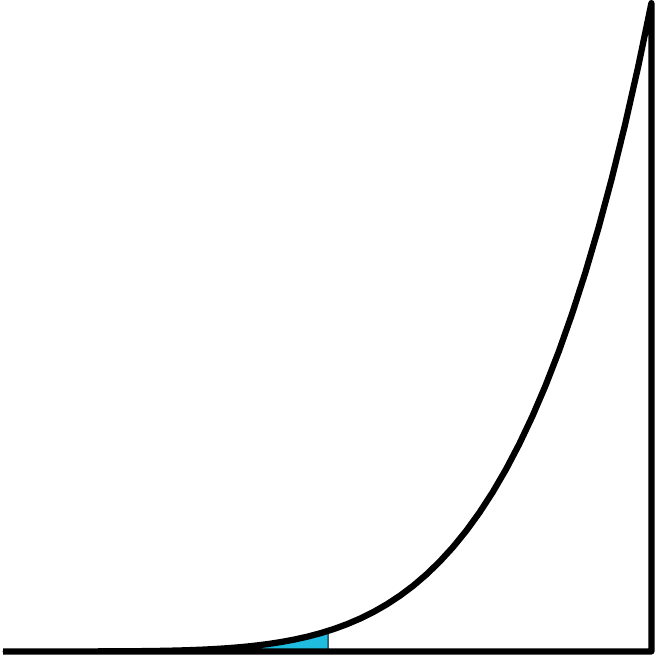}
\subcaption{$\omega_\delta$, $\gamma=5$}
\label{fig:reg:4}
\end{subfigure}
\caption{\emph{Region of regularization.} We here illustrate the region of regularization $\omega_\delta$ for the square-to-cusp parameterization when $\delta=0.001$. On the left, we view the region in the reference domain where the different shades indicate three different $\gamma$. The three illustrations on the right are the corresponding regions in the physical domain. Note that the size of $\hat{\omega}_\delta$ increases quite noticeably with increasing $\gamma$.}
\label{fig:cuspref}
\end{figure}

\paragraph{Model Singular Map.} We continue the derivation in the specific case of the square-to-cusp parameterization detailed in Section~\ref{sec:geometry} and illustrated in Figure~\ref{fig:cuspmap}. Recall the following leading order terms for the eigenvalues and the measure
\begin{align} \label{eq:cuspmap-basic-info}
\lambda_1\sim 1 ,\quad
\lambda_2 \sim \hatx^{2\gamma} ,\quad
|G|^{1/2}\sim \hatx^\gamma
\end{align}
For this parameterization, the set  
$\hat{\omega}_\delta$ is approximately the strip
\begin{equation}
[0,\delta^{1/(2\gamma)}]\times [0,1]
\end{equation}
since we require $\delta^{1/2} \geq \lambda_2^{1/2} \sim (\hatx^{2\gamma})^{1/2} \sim \hatx^\gamma$. This set is illustrated in Figure~\ref{fig:cuspref}.
The square of the remaining factor in \eqref{eq:basic-est} hence scales as
\begin{align}
\bigl\| \lambda_j^{-1/4} |G|^{1/4} \bigr\|_{L^2(\hat{\omega}_\delta)}^2
&=
\int_{\hat{\omega}_\delta}
\lambda_j^{-1/2}  |G|^{1/2}
\sim
\int_{\hat{\omega}_\delta} 1
\sim
\int_{0}^{1}
\int_{0}^{\delta^{1/(2\gamma)}} 1
\sim \delta^{\frac{1}{2\gamma}}
\end{align}

In summary, this gives the consistency error estimate
\begin{align}
\int_{\hatO} \bigl(|G|^{1/2} G^{-1}  \hatnabla \hatu \bigr) \cdot \hatnabla \hatv
-
\int_{\hatO} \bigl(\Reg  \hatnabla \hatu \bigr) \cdot \hatnabla \hatv
&\lesssim 
\underbrace{ \delta^{1/4}\delta^{\frac{1}{4\gamma}} }_{=\delta^{(\gamma+1)/(4\gamma)}}
\left\| \nabla u \right\|_{L^\infty(\omega_\delta)} \tn \hatv \tn_{h,\delta}
\end{align}
For an optimal-order method, we need this to scale as $\sim h^p$, implying that we should choose the regularization parameter $\delta$ such that
\begin{equation} \label{eq:delta-bound}
\boxed{\delta \lesssim h^{4\gamma p/(\gamma + 1)}}
\end{equation}

\begin{rem}[Region of Regularization]
Whereas the bounds for $\delta$ above are devised for optimal order approximation properties, this bound is typically suboptimal for optimal order scaling of the stiffness matrix condition number, cf. \cite{MR4394710} where similar bounds are computed to scale an inconsistent type of stabilization. An optimal bound for the condition number would likely involve a somewhat less aggressive scaling of $\delta$, and the region of regularization $\hat{\omega}_\delta$ using such a bound is reasonably the region where the method will run into numerical issues. 
By the illustration of $\hat{\omega}_\delta$ in Figure~\ref{fig:cuspref}, and taking into account that $\delta$ will scale with $h$, we realize that for smaller $\gamma$ this region will be very small, perhaps only covering only a fraction of each element adjacent to $\hatx=0$. This is likely the reason why quadrature using Gauss points is often sufficient without regularization, as it is unlikely that a Gauss point will be placed inside the region of regularization. For higher $\gamma$, the region grows quite a bit, which increases the likelihood of failure. In addition, our weak enforcement of interface and boundary conditions increases the likelihood that quadrature points are placed within the region of regularization, which means regularization is even more critical for robustness.
\end{rem}

\begin{figure}
\centering
\begin{subfigure}[t]{0.23\linewidth}\centering	
\includegraphics[width=0.95\linewidth]{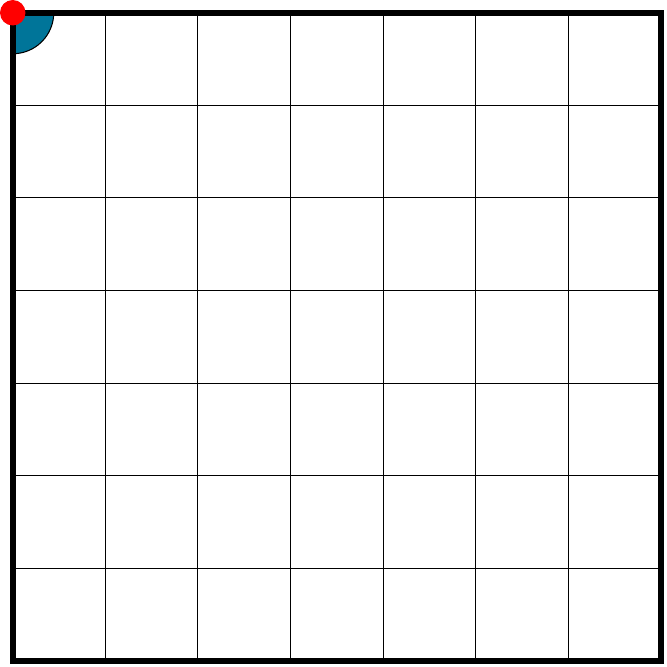}
\subcaption{$\widehat{\Omega}$}
\label{fig:bilinear:1}
\end{subfigure}
	\quad
\begin{subfigure}[t]{0.23\linewidth}\centering
\includegraphics[width=0.95\linewidth]{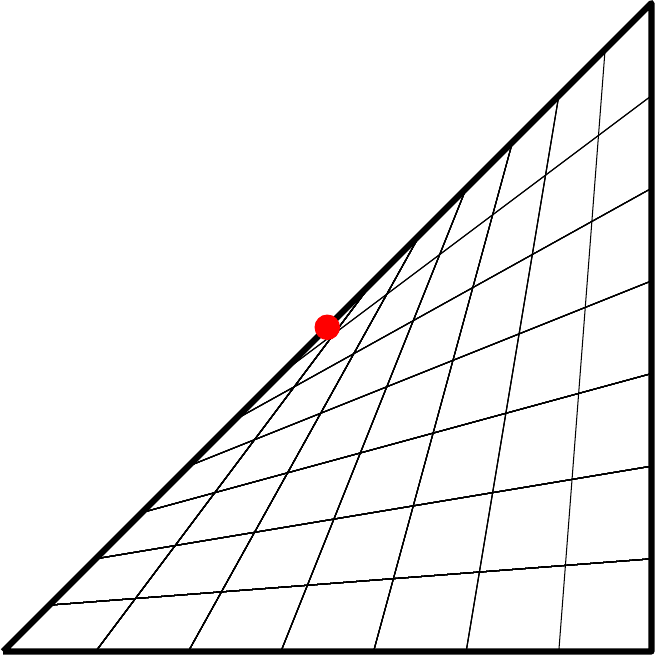}
\subcaption{$\Omega$}
\label{fig:billinear:2}
\end{subfigure}
\caption{\emph{Collapsed bilinear map.} Illustration of a singular parameterization where the singularity stems from the map derivatives becoming colinear at a certain point. The region of regularization $\hat{\omega}_\delta$ is indicated for $\delta=0.001$.}
\label{fig:bilinear}
\end{figure}

\paragraph{Colinear Map Derivatives.}
Another common type of singular parameterization is when a corner in the reference domain is mapped onto a point on the physical boundary where the boundary is smooth. This leads to colinear rows in the map derivative when approaching the corner with a singular metric, cf. \cite{MR2576766,LU20092391}. A simple example is the bilinear map illustrated in Figure~\ref{fig:bilinear} where the upper left corner of the unit square is mapped onto the midpoint of the square diagonal.
This map is expressed
\begin{equation}
F
\left(
\left[
\begin{matrix}
\hatx
\\
\haty
\end{matrix}
\right]
\right)
=
\left[
\begin{matrix}
 \hatx + \frac{(1-\hatx)\haty}{2}
\\
\hatx\haty + \frac{(1-\hatx)\haty}{2}
\end{matrix}
\right]
\end{equation}
and its Riemannian metric tensor is
\begin{align}
G =
\frac{1}{2}
\begin{bmatrix}
\haty^2-2\haty+2 & \hatx\haty - \hatx + 1 \\
\hatx\haty - \hatx + 1 &         \hatx^2 + 1
\end{bmatrix}
\end{align}
which becomes singular as $(\hatx,\haty)$ approaches the corner $(0,1)$.
By similar arguments as in the square-to-cusp case, we
find that the leading terms for the eigenvalues to $G$ are
\begin{equation}
\lambda_1 \sim 1,\qquad \lambda_2 \sim \frac{1}{4}\hatr^2
\end{equation}
where $\hatr = \sqrt{\hatx^2 + (\haty-1)^2}$ is the radial distance from the singular point $(0,1)$. This implies a region of regularization
\begin{align}
\hat{\omega}_\delta = \bigl\{ (\hatx,\haty) \in \hatO \, | \, \hatr < 2\delta^{1/2} \bigr\}
\end{align}
Computing the square of the remaining factor in \eqref{eq:basic-est} gives us
\begin{align}
\bigl\| \lambda_i^{-1/4} |G|^{1/4} \bigr\|_{L^2(\hat{\omega}_\delta)}^2
=
\int_{\hat{\omega}_\delta} \lambda_2^{-1/2} |G|^{1/2}
\sim
\int_{\hat{\omega}_\delta} 1
\sim
\int_{\hatr=0}^{\hatr=2\delta^{1/2}} 1
\sim
\int_{0}^{2\delta^{1/2}} \hatr
\sim \delta
\end{align}
where we, in the final integral, make a change of variables to polar coordinates.
Accounting for the factor $\delta^{1/4}$ in \eqref{eq:basic-est}, this implies that the regularization parameter $\delta$ should be chosen such that
\begin{equation} \label{eq:delta-colinear-bound}
\boxed{\delta \lesssim h^{4p/3}}
\end{equation}
for an optimal order method, which we note is less aggressive than the bound \eqref{eq:delta-bound} for the square-to-cusp map.

\section{Numerical Results} \label{sec:numerics}
In this section, we demonstrate the robustness and performance of our method on multipatch domains, which include patches described using the singular square-to-cusp parameterization detailed in Section~\ref{sec:model-singular-map}.

\subsection{Implementation}

A 2D version of the method \eqref{eq:the-method} was implemented in {\sc Matlab}. This was largely based on our prior implementations in \cite{MR3682761, MR3709202} with the regularization as the main addition, and we refer to those papers for more details on implementational aspects.

\paragraph{Geometry Description.}
The geometry of each patch $\Omega_i$ is described as a mapping $F_i:\hatO_i\to\Omega_i$ and trim curves defining the boundary of the reference domain $\partial\hatO_i$. In our implementation, we approximate $\partial\hatO_i$ by a polygonal boundary of a resolution that can be chosen independently from the approximation space. If the trim curves are curved, this leads to geometric errors that may be controlled by increasing the resolution of the polygonal boundary. However, in all numerical examples below, the reference domain is the unit square, which means that $\partial\hatO_i$ is represented exactly using this geometry description.

\paragraph{Patch Coupling.}
The coupling of the numerical solution over the patch interfaces is in the method \eqref{eq:the-method} largely hidden in the average $\langle \hatu \rangle$ defined by \eqref{eq:average}.
In our setting, the patch connectivity is known from the geometric description, where each interface is described via mappings from $[0,1]$ to the corresponding points on the respective patches joined at the interface.

\paragraph{Approximation Space.}
We place each patch reference domain $\hatO_i$ on top of a structured quadrilateral background grid with element size $h$ and let the computational mesh $\mcK_{h,i}$ be constituted of all quadrilaterals in the grid with a non-zero intersection with $\hatO_i$.
On each mesh, we define an approximation $\hatV_{h,i}$ space using full regularity tensor product B-splines of degree $p=1,2,3$, which gives the patchwise regularity $C^{p-1}(\hatO_i) \subset H^p(\hatO_i)$.

\paragraph{Quadrature.}
If the mesh $\mcK_{h,i}$ is \emph{matching}, i.e., that the mesh is a partition of $\hatO_i$, we use tensor products of 1D Gauss rules to integrate the quadrilateral elements. If the mesh is \emph{cut}, i.e., that the mesh is not matching and only covers $\hatO_i$, the cut elements that are clipped by the reference domain boundary $\partial\hatO_i$ are treated differently. For these elements, their intersection with $\hatO_i$ is triangulated, whereafter a Gauss rule on each triangle is applied. For the boundary integrals, 1D Gauss rules are used. By this construction, a quadrature point should never be placed precisely on the singular point, but they may be close.

\paragraph{Parameter Choices.}
The Nitsche penalty parameter is chosen as $\beta=25 p^2$ and the ghost penalty stabilization parameter as $\eta = 0.01$. Unless otherwise stated, ghost penalty stabilization is only applied when using trimmed approximation spaces. The weights in the average \eqref{eq:average} are set to $\kappa_{ij}=0.5$, giving equal weight to the neighboring patches.

\subsection{Model Problem} \label{sec:model-problem}
\begin{figure}
	\centering
	\begin{subfigure}[t]{0.23\linewidth}\centering	
		\includegraphics[width=0.95\linewidth]{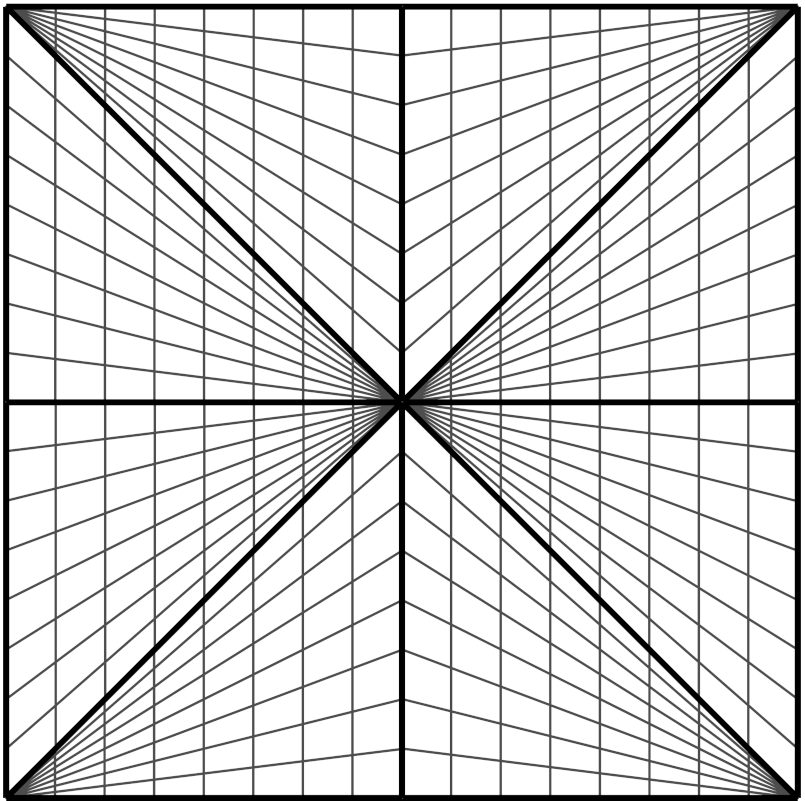}
		\subcaption{$\gamma=1$, matching}
		\label{fig:meshes:1}
	\end{subfigure}
	\begin{subfigure}[t]{0.23\linewidth}\centering
		\includegraphics[width=0.95\linewidth]{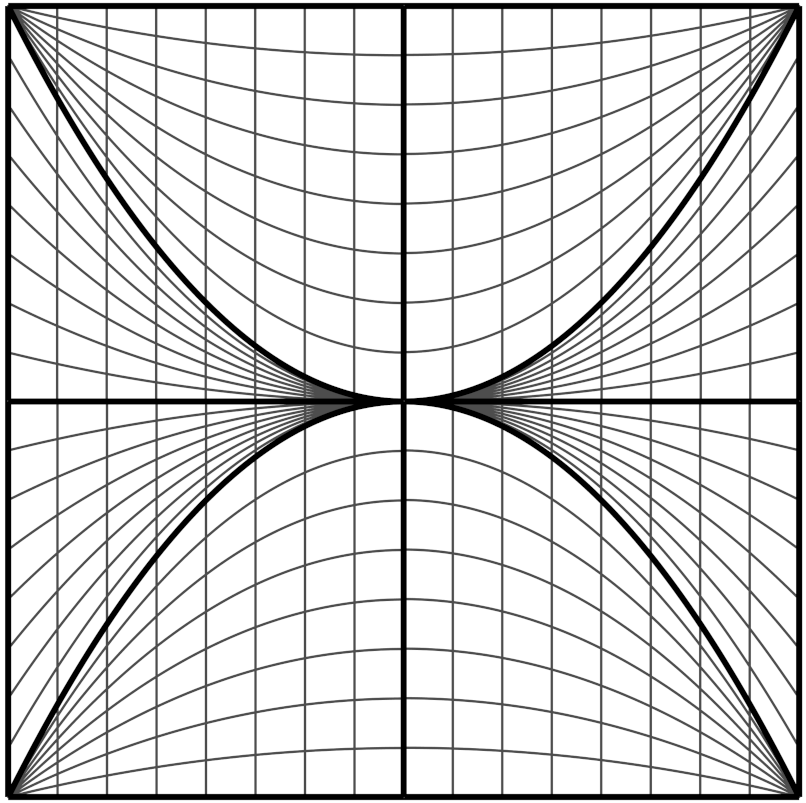}
		\subcaption{$\gamma=2$, matching}
		\label{fig:meshes:2}
	\end{subfigure}
	\begin{subfigure}[t]{0.23\linewidth}\centering
		\includegraphics[width=0.95\linewidth]{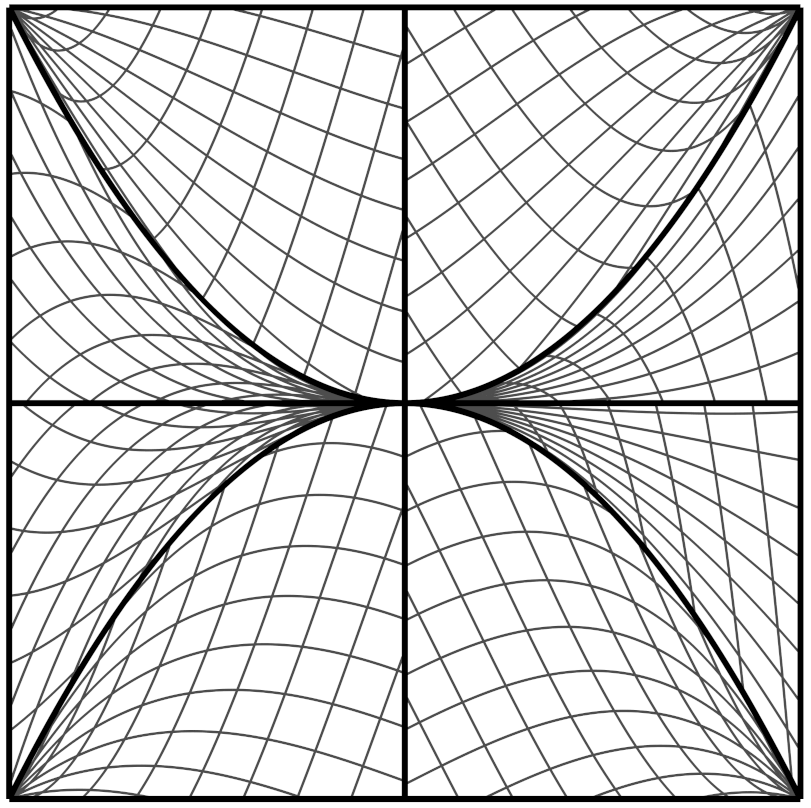}
		\subcaption{$\gamma=2$, trimmed}
		\label{fig:meshes:3}
	\end{subfigure}
	\begin{subfigure}[t]{0.23\linewidth}\centering	
			\includegraphics[width=0.95\linewidth]{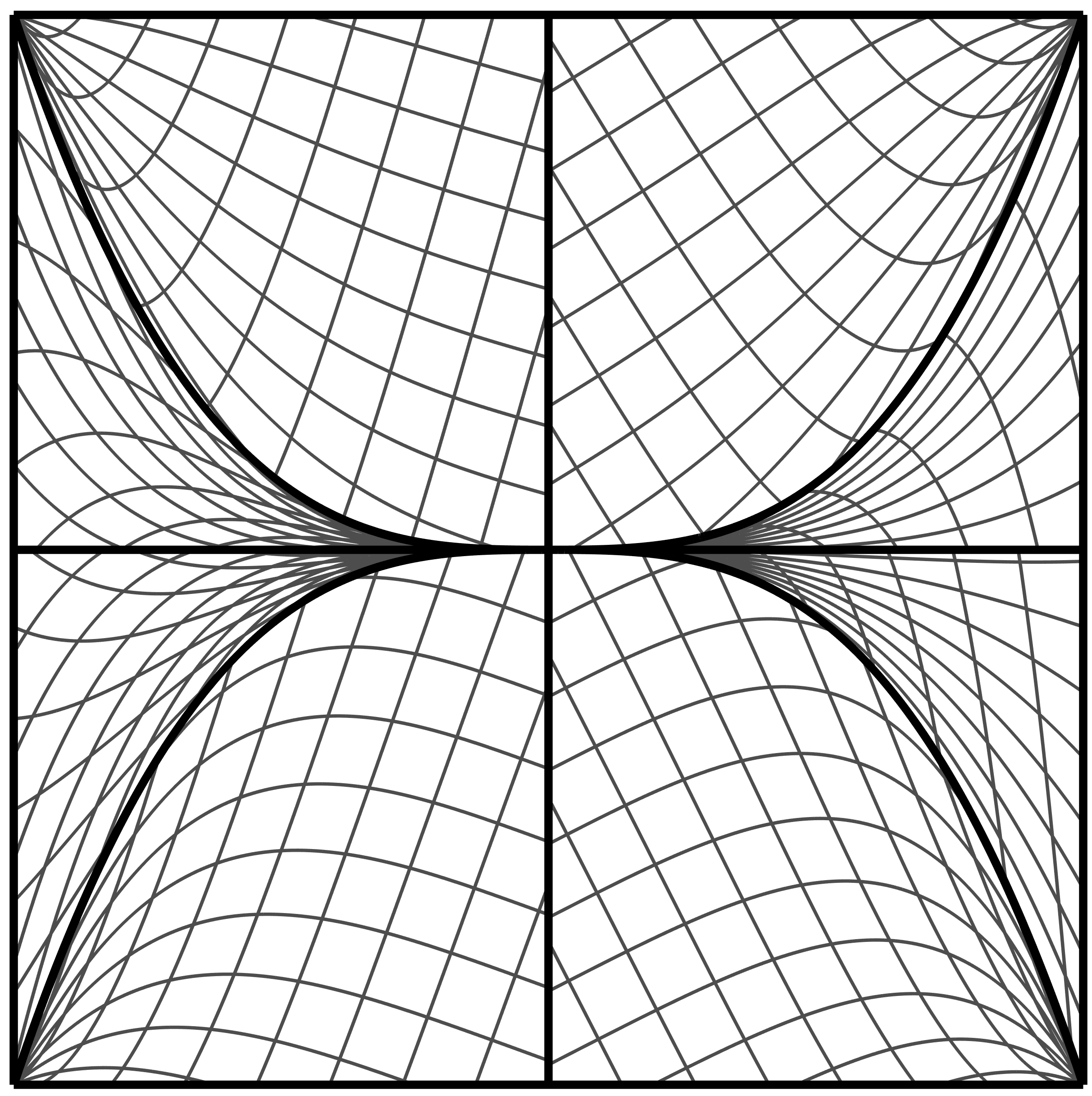}
			\subcaption{$\gamma=3$, trimmed}
			\label{fig:nsmooth}
		\end{subfigure}
\caption{\emph{Model multipatch domain with meshes.} Four different realizations of the model multipatch domain ($\gamma=1,2,3$). The two realizations on the left are equipped with quadrilateral meshes that match the patch reference domains, while the two realizations on the right are equipped with trimmed (cut) meshes in the patch reference domains constructed by different rotations of the background grid from which the meshes are extracted.}
\label{fig:meshes}
\end{figure}

As a model problem on a multipatch domain $\Omega$ with singular parameterizations, we consider a partition of the square $[-1,1]^2$ into 8 patches $\Omega_i$, where each patch is defined as a map $F_i : \widehat{\Omega}_i=[0,1]^d \rightarrow \Omega_i$. 
The maps associated with the patches are
\begin{align}
F_1(\hatx,\haty) &= \big(  {\hatx},\hatx^\gamma \haty \big)  &	F_2(\hatx,\haty) &= \big( \hatx,(1-\haty) \hatx^\gamma+\haty \big)
\\
F_3(\hatx,\haty) &= \big(  \hatx-1,\haty(1-\hatx)^\gamma \big)  &	F_4(\hatx,\haty) &= \big(  -\hatx,1-\haty(1-\hatx^\gamma) \big)
\\
F_5(\hatx,\haty) &= \big(  \hatx,(1-\haty)\hatx^\gamma)  &	F_6(\hatx,\haty) &= \big(  \hatx,\haty(1-\hatx^\gamma)-1 \big)
\\
F_7(\hatx,\haty) &= \big(  -\hatx,-\haty\hatx^\gamma \big)  &	F_8(\hatx,\haty) &= \big(  -\hatx,(1-\haty)(1-\hatx^\gamma)-1 \big)
\end{align}
where we note that several maps have the structure of the square-to-cusp model geometry detailed in Section~\ref{sec:model-singular-map} and that we hence may use the parameter $\gamma$ to vary the aggressiveness of the cusp. Realizations of this multipatch domain for different $\gamma$ are presented in Figure~\ref{fig:meshes}, where we have equipped the unit square reference domains with a structured grid, which is mapped onto the physical domain. This grid is, in some cases, aligned with the unit square reference domain, giving a matching mesh on each patch, and in some cases, the grid is rotated, creating trimmed meshes where quadrilaterals along the reference domain boundary are cut. While this multipatch domain might seem very simple, it very much pinpoints the problematics of singular maps in trimmed multipatch methods.

\begin{figure}
\centering	
\begin{subfigure}[t]{0.32\linewidth}\centering
\includegraphics[width=0.9\linewidth]{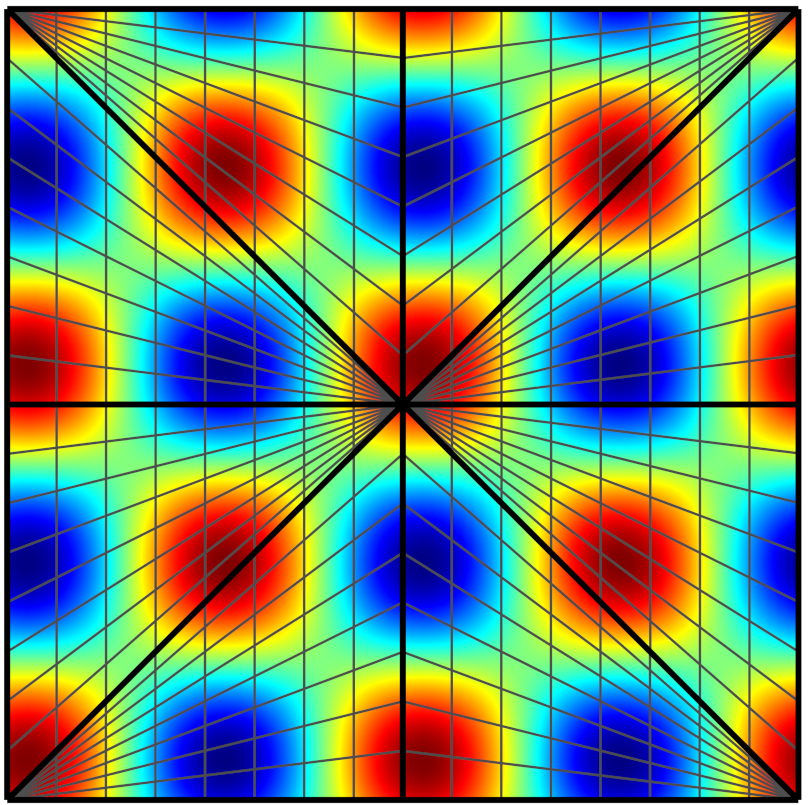}
\subcaption{$\gamma=1$, matching}
\label{fig:numerical-solution-a}
\end{subfigure}
\begin{subfigure}[t]{0.32\linewidth}\centering
\includegraphics[width=0.9\linewidth]{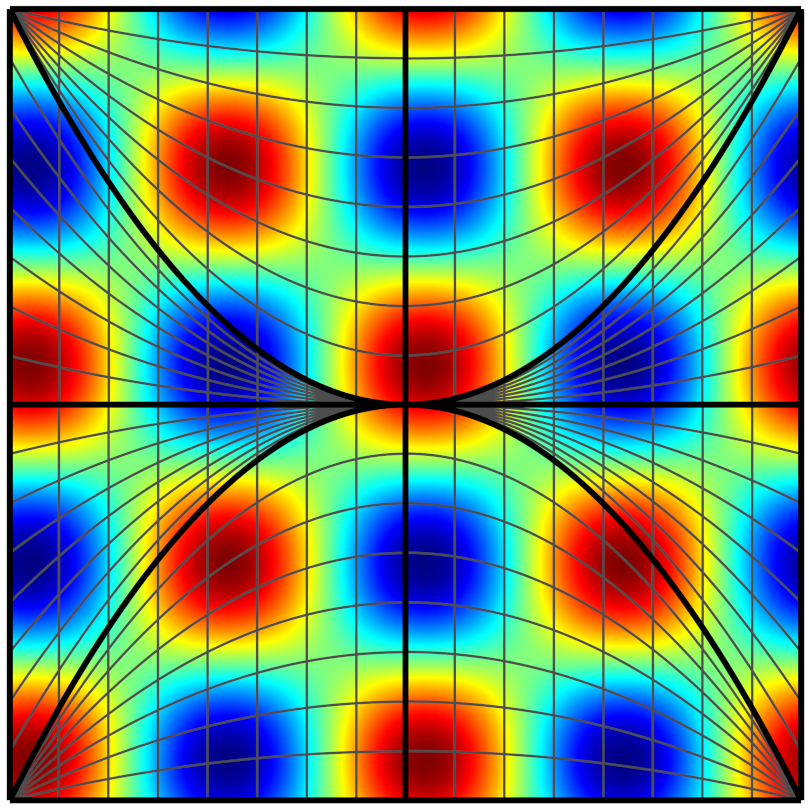}
\subcaption{$\gamma=2$, matching}
\label{fig:numerical-solution-b}
\end{subfigure}
\begin{subfigure}[t]{0.32\linewidth}\centering
\includegraphics[width=0.9\linewidth]{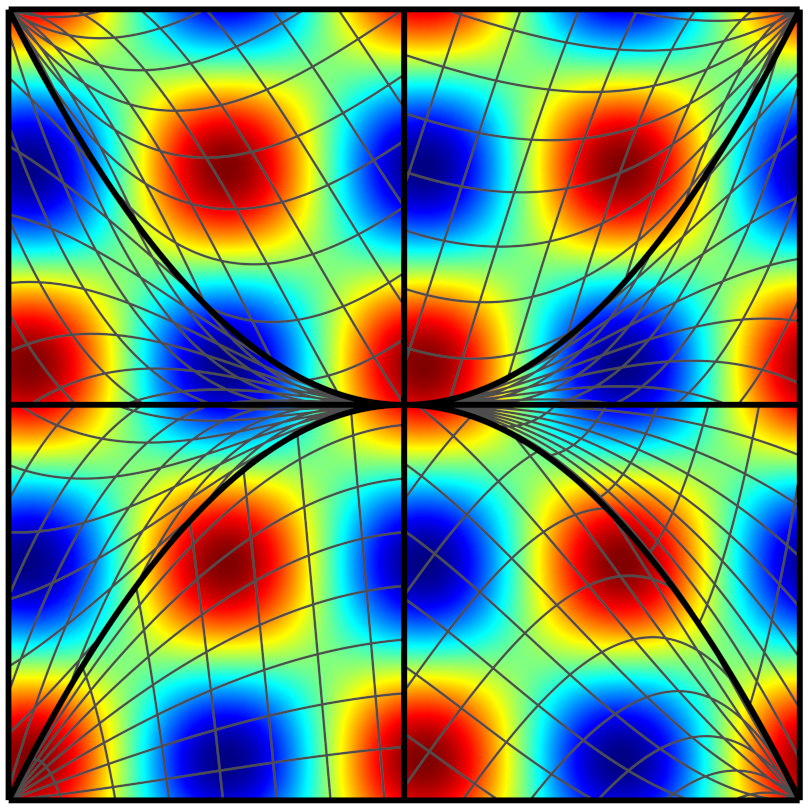}
\subcaption{$\gamma=2$, trimmed}
\label{fig:numerical-solution-c}
\end{subfigure}

\caption{\emph{Sample numerical solutions.} Numerical solutions to the model problem with cusp aggressiveness $\gamma=1,2$ using approximation spaces of maximum regularity B-splines ($p=2$).}
\label{fig:numerical-solution}
\end{figure}

We manufacture a problem with known analytical solution on $\Omega$ by making the ansatz $u =  \sin(2\pi (x-0.3)) \cos(2\pi (y+0.4))$, and defining the patchwise solution as the restriction $u_i = u|_{\Omega_i}$, from which we derive the data $f_i$ and $g_i$. In Figure~\ref{fig:numerical-solution}, we present sample numerical solutions for this model with $\gamma=1,2$ using matching respectively trimmed (cut) approximation spaces.

\subsection{Experiments}
We here present the results from our numerical experiments.
The norms $\| u - u_h \|$ and $\| \nabla(u - u_h) \|$ we compute in these experiments are the standard $L^2(\Omega)$-norm and the regularized $H^1(\Omega)$-seminorm, i.e.,
\begin{align}
\| u - u_h \|^2 &= \sum_{i\in\mcI}\int_{\hatO_i} (u_i\circ F_i - \hatu_{h,i})^2 |G_i|^{1/2}
\\
\| \nabla ( u - u_h ) \|^2 &= \sum_{i\in\mcI}\int_{\hatO_i} \bigl(\Regi\hatnabla (u_i\circ F_i - \hatu_{h,i})\bigr) \cdot \hatnabla (u_i\circ F_i - \hatu_{h,i})
\end{align}

\begin{figure}
\centering
\begin{subfigure}[t]{0.32\linewidth}\centering
\includegraphics[width=1\linewidth]{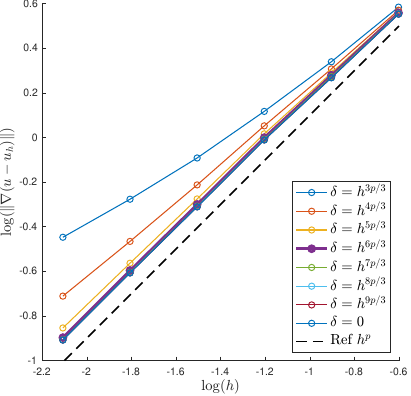}
\subcaption{$\gamma=1$, $p=1$}
\label{fig:delta-scaling-p1-gamma1-H1}
\end{subfigure}
\begin{subfigure}[t]{0.32\linewidth}\centering
\includegraphics[width=1\linewidth]{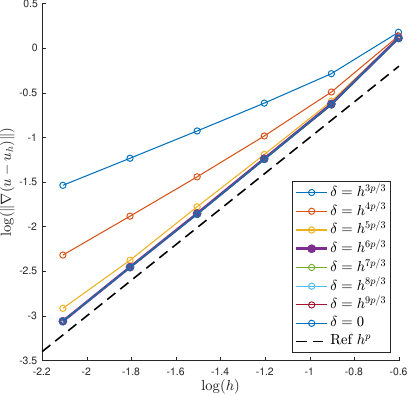}
\subcaption{$\gamma=1$, $p=2$}
\label{fig:delta-scaling-p2-gamma1-H1}
\end{subfigure}
\begin{subfigure}[t]{0.32\linewidth}\centering
\includegraphics[width=1\linewidth]{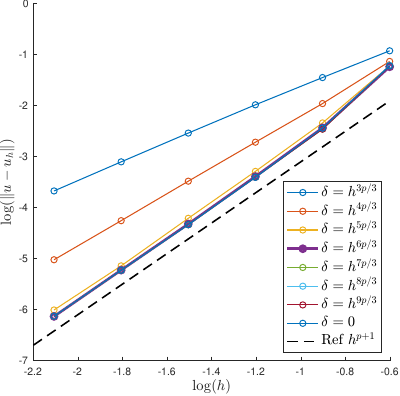}
\subcaption{$\gamma=1$, $p=2$}
\label{fig:delta-scaling-p2-gamma1-L2}
\end{subfigure}
\\[1ex]

\begin{subfigure}[t]{0.32\linewidth}\centering
\includegraphics[width=1\linewidth]{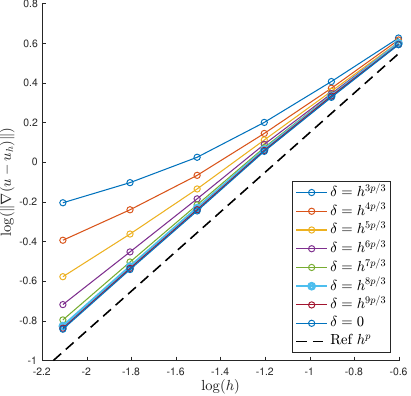}
\subcaption{$\gamma=2$, $p=1$}
\label{fig:delta-scaling-p1-gamma2-H1}
\end{subfigure}
\begin{subfigure}[t]{0.32\linewidth}\centering
\includegraphics[width=1\linewidth]{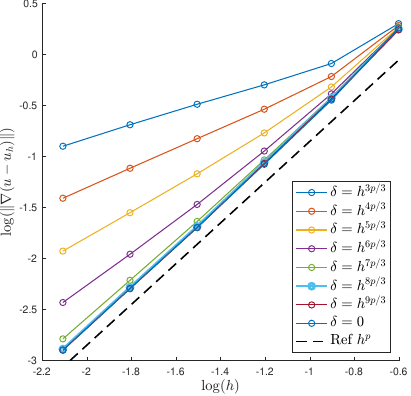}
\subcaption{$\gamma=2$, $p=2$}
\label{fig:delta-scaling-p2-gamma2-H1}
\end{subfigure}
\begin{subfigure}[t]{0.32\linewidth}\centering
\includegraphics[width=1\linewidth]{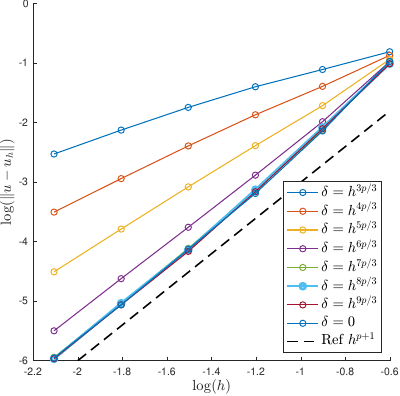}
\subcaption{$\gamma=2$, $p=2$}
\label{fig:delta-scaling-p2-gamma2-L2}
\end{subfigure}
\\[1ex]

\begin{subfigure}[t]{0.32\linewidth}\centering
\includegraphics[width=1\linewidth]{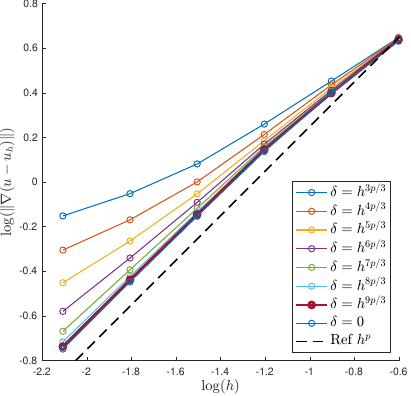}
\subcaption{$\gamma=3$, $p=1$}
\label{fig:delta-scaling-p2-gamma5-H1}
\end{subfigure}
\begin{subfigure}[t]{0.32\linewidth}\centering
\includegraphics[width=1\linewidth]{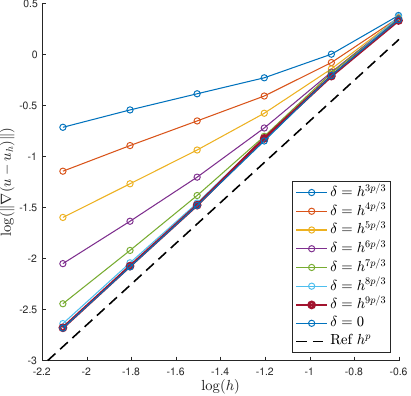}
\subcaption{$\gamma=3$, $p=2$}
\label{fig:delta-scaling-p2-gamma5-H1}
\end{subfigure}
\begin{subfigure}[t]{0.32\linewidth}\centering
\includegraphics[width=1\linewidth]{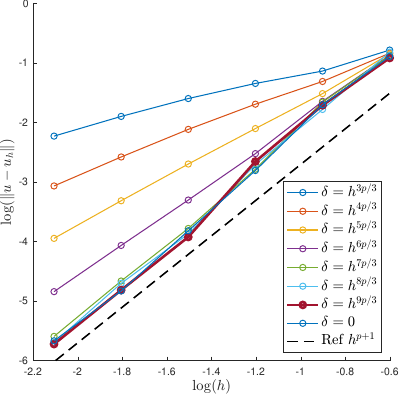}
\subcaption{$\gamma=3$, $p=2$}
\label{fig:delta-scaling-p2-gamma5-L2}
\end{subfigure}

\caption{\emph{Scaling of the regularization parameter $\delta$.} Numerical study on which scaling of the regularization parameter $\delta$ is needed for optimal order convergence for the multipatch model problem in $H^1$-seminorm and $L^2$-norm using cusp aggressiveness $\gamma\in\{1,2,3\}$. We use matching meshes in each reference patch equipped with maximal regularity B-splines of orders $p=1,2$. The thicker line in the plots indicates the scaling implied by our consistency estimate, $\delta \lesssim h^{4\gamma p/(\gamma + 1)}$, and the results suggest this scaling is sufficient for optimal order approximation.
}
\label{fig:delta-scaling}
\end{figure}

\paragraph{Regularization Parameter Scaling.}
To assess the validity of our proposed bound \eqref{eq:delta-bound}, where we state that the regularization parameter needs to scale as $\delta \lesssim h^{4\gamma p/(\gamma + 1)}$ for an optimal order method, we test the convergence behavior of the regularized method for a range of $\delta$-scalings using different cusp aggressiveness $\gamma$ and approximation space orders $p$.
Results from this experiment are collected in Figure~\ref{fig:delta-scaling}. These results indicate that our proposed bound is sharp, i.e., that this scaling is sufficient for an optimal order method.

\paragraph{Convergence and Robustness.}
We continue by studying convergence and stiffness matrix condition number behavior of the regularized method using $\delta = h^{4\gamma p/(\gamma + 1)}$, satisfying the bound \eqref{eq:delta-bound}, for the model problem with cusp aggressiveness $\gamma=2$.

\begin{figure}
\centering	
\begin{subfigure}[t]{0.32\linewidth}\centering
	\includegraphics[width=1\linewidth]{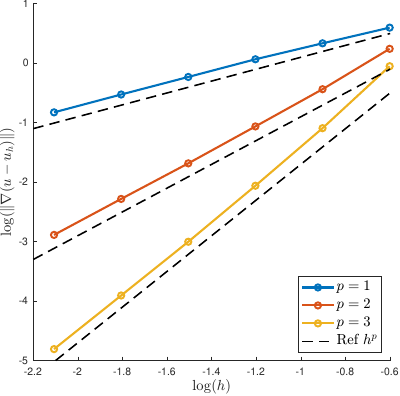}
\subcaption{$\gamma=2$, matching}
\label{fig:convergence-matching-H1}
\end{subfigure}
\begin{subfigure}[t]{0.32\linewidth}\centering
	\includegraphics[width=1\linewidth]{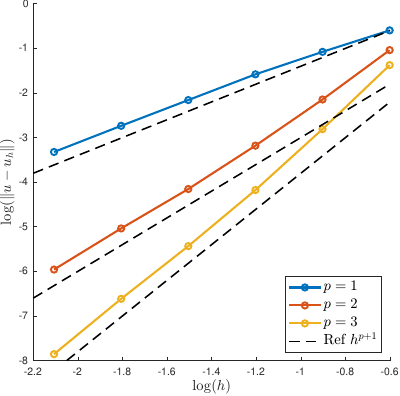}
\subcaption{$\gamma=2$, matching}
\label{fig:convergence-matching-L2}
\end{subfigure}
\\[1ex]	

\begin{subfigure}[t]{0.32\linewidth}\centering
	\includegraphics[width=1\linewidth]{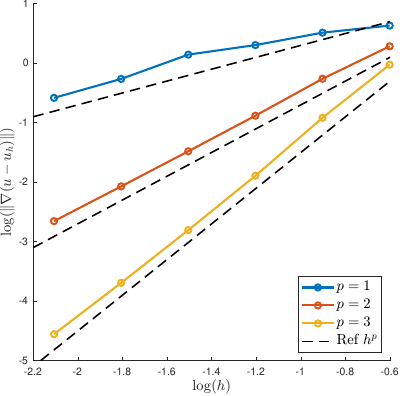}
\subcaption{$\gamma=2$, trimmed}
\label{fig:convergence-cut-H1}
\end{subfigure}
\begin{subfigure}[t]{0.32\linewidth}\centering
	\includegraphics[width=1\linewidth]{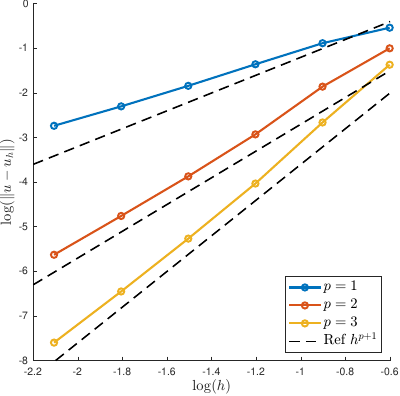}
\subcaption{$\gamma=2$, trimmed}
\label{fig:convergence-cut-L2}
\end{subfigure}

\caption{\emph{Convergence.} We here study the convergence behavior of the regularized method ($\delta = h^{4\gamma p/(\gamma + 1)}$) applied to the model problem for moderately aggressive cusps using both matching and trimmed (cut) approximation spaces. We note optimal order convergence in $H^1$-seminorm and $L^2$-norm for all studied cases.}
\label{fig:convergence}
\end{figure}

In Figure~\ref{fig:convergence}, we present convergence results using approximation spaces of order $p=1,2,3$ on both matching and trimmed reference domain meshes, and we note optimal order convergence in both $L^2(\Omega)$-norm and in $H^1(\Omega)$-seminorm.

\begin{figure}
	\centering	
	\begin{subfigure}[t]{0.32\linewidth}\centering
		\includegraphics[width=1\linewidth]{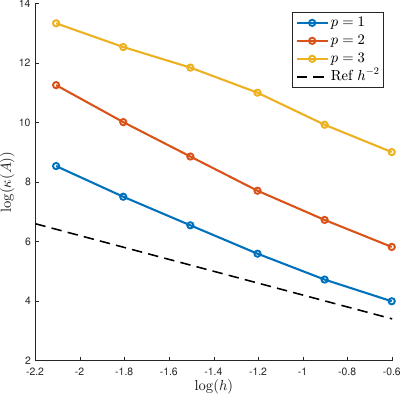}
	\subcaption{$\gamma=2$, matching}
	\label{fig:condition-number:a}
	\end{subfigure}
	\begin{subfigure}[t]{0.32\linewidth}\centering
		\includegraphics[width=1\linewidth]{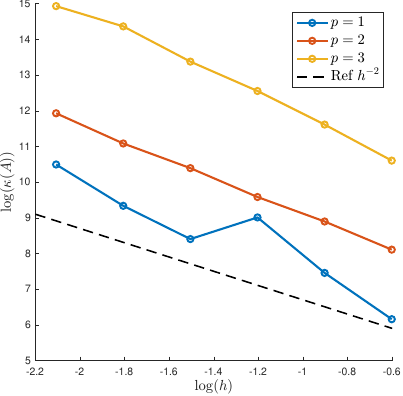}
	\subcaption{$\gamma=2$, trimmed}
	\label{fig:condition-number:b}
	\end{subfigure}
	\begin{subfigure}[t]{0.32\linewidth}\centering
	\includegraphics[width=0.95\linewidth]{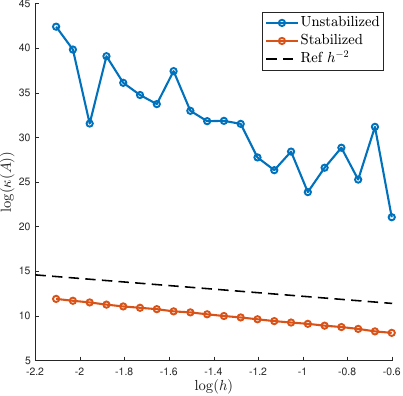}
	\subcaption{$\gamma=2$, $p=2$, trimmed}
	\label{fig:condition-number:c}
\end{subfigure}
\caption{\emph{Condition number scaling.} We here study how the stiffness matrix condition number of the regularized method ($\delta = h^{4\gamma p/(\gamma + 1)}$) scales with the mesh size $h$ for the model problem ($\gamma=2$) using matching and trimmed approximation spaces. We note that the scaling seems to be slightly worse than $\sim h^{-2}$ in this range of mesh sizes. In (c), we illustrate the importance of including the ghost penalty stabilization \eqref{eq:ghost-penalty} when using a trimmed (cut) approximation space.}
\label{fig:condition-number}
\end{figure}

In Figure~\ref{fig:condition-number}, we study how the stiffness matrix condition number scales with $h$ on both matching and trimmed reference domain meshes. The scaling within this range of mesh sizes is slightly worse than the desired $h^{-2}$-scaling, which is reasonable since the regularization parameter $\delta$ was chosen for optimal order accuracy rather than optimal condition number scaling. We assume that for sufficiently small $h$, the method will approach a $h^{-2}$ scaling since the effect of the regularization will diminish.

To illustrate the necessity of including ghost penalty stabilization \eqref{eq:ghost-penalty} when using trimmed approximation spaces, we in Figure~\ref{fig:condition-number:c} include a study of the condition number scaling for trimmed $p=2$ approximation spaces, with and without stabilization. Depending on the cut situations caused by the trimming, without stabilization the condition number may grow arbitrarily large and we can also have loss of coercivity.

\paragraph{Extreme Singularities.}
Next, we investigate the behavior of the method in the model problem with higher values of $\gamma$, yielding more extreme singularities.
Matching approximation spaces of order $p=2$ are used.
We consider three variants of the method; the regularized method with $\delta = h^{4\gamma p/(\gamma + 1)}$ fulfilling our bound \eqref{eq:delta-bound}; the regularized method with $\delta=0$, which is actually not regularized, but computes $|G_i|^{1/2} G_i^{-1}$ more robustly; and a naive implementation of the unregularized method, computing $|G_i|^{1/2} G_i^{-1}$ as a product of two terms.

\begin{figure}
\centering
\begin{subfigure}[t]{0.32\linewidth}\centering
	\includegraphics[width=0.95\linewidth]{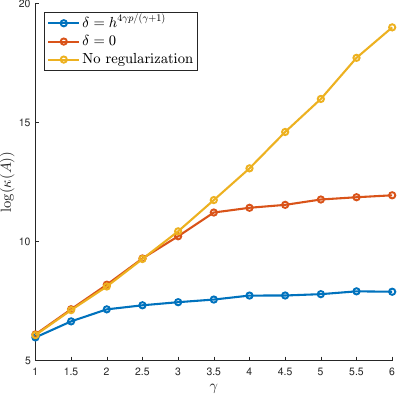}
	\subcaption{Condition number, $\gamma\in[1,6]$}
	\label{fig:extreme-cusps:a}
\end{subfigure}
\begin{subfigure}[t]{0.32\linewidth}\centering
	\includegraphics[width=0.95\linewidth]{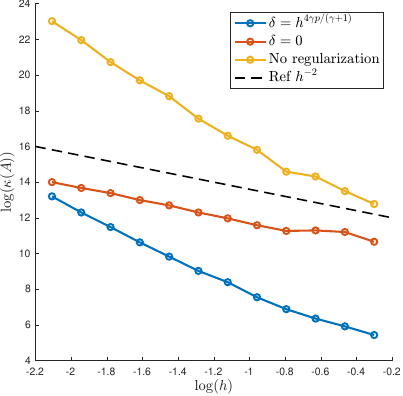}
	\subcaption{Condition number, $\gamma=5$}
	\label{fig:extreme-cusps:b}
\end{subfigure}
\\[1ex]

\begin{subfigure}[t]{0.32\linewidth}\centering
	\includegraphics[width=0.95\linewidth]{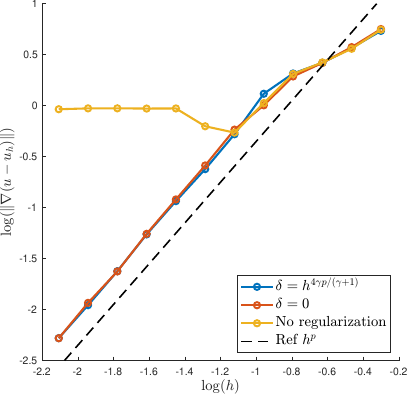}
	\subcaption{$H^1$ error, $\gamma=5$}
	\label{fig:extreme-cusps:c}
\end{subfigure}
\begin{subfigure}[t]{0.32\linewidth}\centering
	\includegraphics[width=0.95\linewidth]{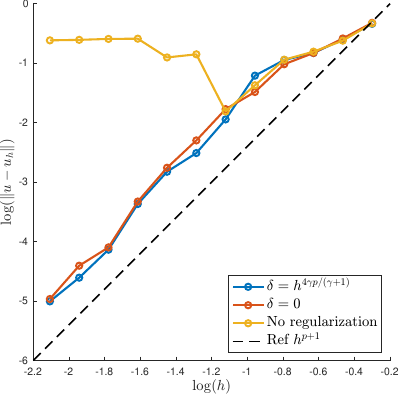}
	\subcaption{$L^2$ error, $\gamma=5$}
	\label{fig:extreme-cusps:d}
\end{subfigure}
	
\caption{\emph{Extreme singularities.} We here study the convergence and stiffness matrix condition number behavior for the model problem with more aggressive cusps. Approximation spaces with matching meshes and $p=2$ are used, and we consider the regularized method with $\delta = h^{4\gamma p/(\gamma + 1)}$ according to our bound \eqref{eq:delta-bound} and $\delta=0$. In addition, we include a naive implementation of the unregularized method, computing $|G_i|^{1/2} G_i^{-1}$ as a product of two terms. In (a), we consider the condition number for a fixed mesh size $h=0.1$ when gradually increasing cusp aggressiveness $\gamma$. In (b), we investigate how the condition number scales with $h$ when $\gamma=5$, and in (c)--(d), we look at the convergence with respect to $h$ when $\gamma=5$.}
\label{fig:extreme-cusps}
\end{figure}

In Figure~\ref{fig:extreme-cusps:a}, we study the behavior of the stiffness matrix condition number on a fixed mesh size $h=0.1$ when gradually increasing the cusp aggressiveness from $\gamma=1$ to $\gamma=6$. We note that the naive method without regularization behaves very poorly, with a rapidly increasing condition number, while the regularized method is better behaved with slowly increasing condition numbers. The $\delta=0$ variant also has slowly increasing condition numbers, albeit on a higher level.

In Figure~\ref{fig:extreme-cusps:b}, we consider the condition number scaling of the three method variants when the cusp aggressiveness is $\gamma=5$. Again, the naive method without regularization behaves very poorly with an estimated scaling of roughly $h^{-6}$. The $\delta=0$ variant exhibits the desired $h^{-2}$-scaling, while the regularized method with $\delta = h^{4\gamma p/(\gamma + 1)}$ feature even lower condition numbers but with a worse scaling. Most likely, the regularized method will approach the $\delta=0$ variant with smaller $h$ since $\delta = h^{4\gamma p/(\gamma + 1)}$ becomes very small for smaller $h$, and that effect should be even more rapid with larger $\gamma$.

In figures~\ref{fig:extreme-cusps:c}--\ref{fig:extreme-cusps:d}, we study the convergence of the three variants when $\gamma=5$. We note optimal order convergence for the regularized method and the $\delta=0$ variant, while the naive method without regularization fails to converge for sufficiently small $h$.

\begin{figure}
	\centering	
	\begin{subfigure}[t]{0.32\linewidth}\centering
		\includegraphics[width=0.95\linewidth]{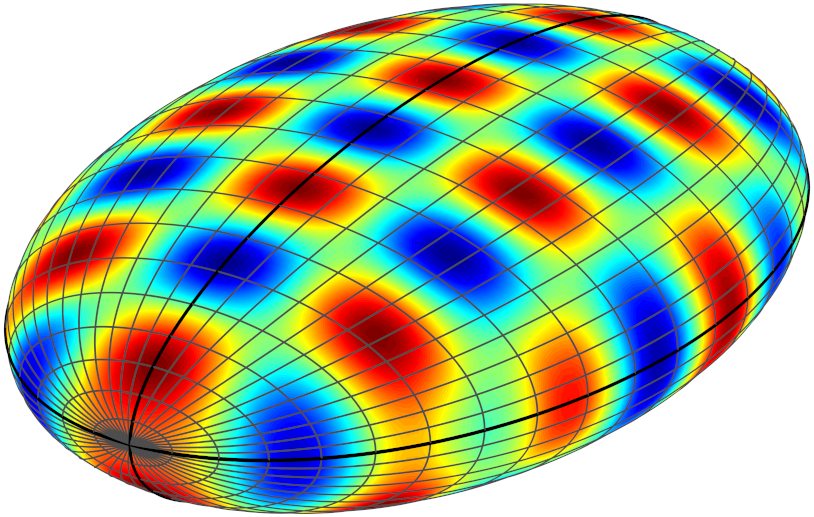}
		\subcaption{Solution}
		\label{fig:ellips:a}
	\end{subfigure}
	\begin{subfigure}[t]{0.32\linewidth}\centering
		\includegraphics[width=0.95\linewidth]{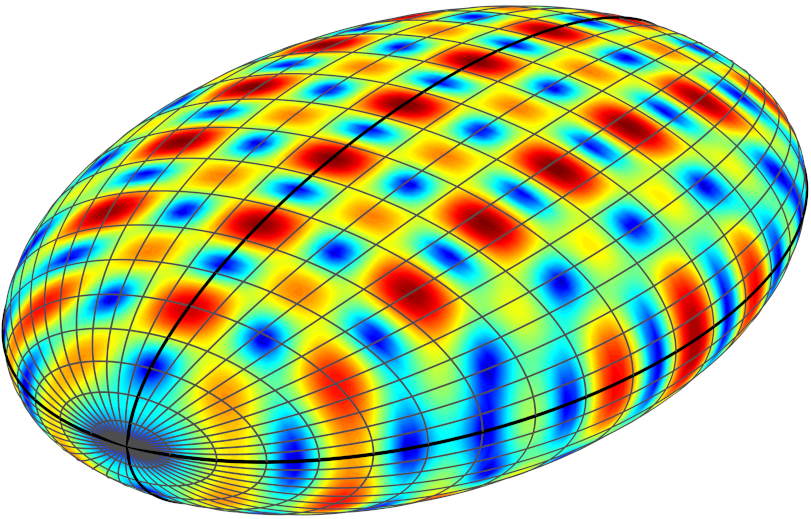}
		\subcaption{Gradient magnitude}
		\label{fig:ellips:b}
	\end{subfigure}
	\\[2ex]

	\begin{subfigure}[t]{0.32\linewidth}\centering
		\includegraphics[width=0.95\linewidth]{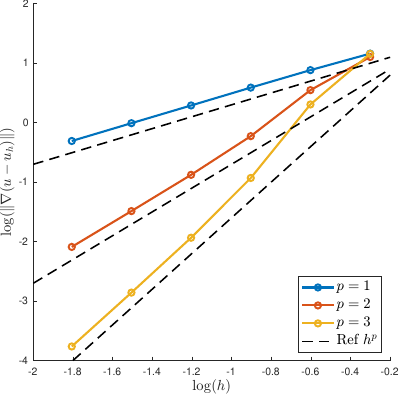}
		\subcaption{$H^1$ error}
		\label{fig:ellips:c}
	\end{subfigure}
	\begin{subfigure}[t]{0.32\linewidth}\centering
		\includegraphics[width=0.95\linewidth]{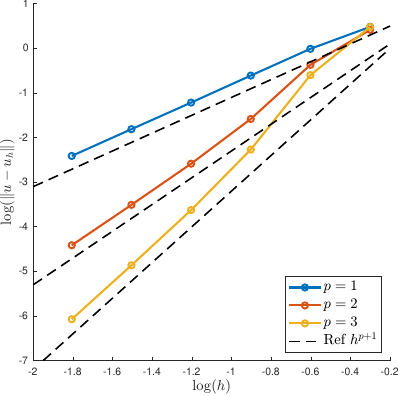}
		\subcaption{$L^2$ error}
		\label{fig:ellips:d}
	\end{subfigure}

\caption{\emph{Numerical solution on a surface.} In (a)--(b), we present the numerical solution of a Laplace-Beltrami problem posed on a multipatch surface, where each of the four patches is described using a singular parameterization. Each patch has two singular points in this case, which correspond to two opposing sides in the unit square reference domains.
In (c)--(d), we present convergence results indicating optimal order convergence.}
\label{fig:ellips}
\end{figure}

\paragraph{Example on a Multipatch Surface.}
To illustrate that the Riemannian setting makes the method readily applicable to problems on domains of positive codimension, we here include a problem on a multipatch surface embedded in $\IR^3$.
The differential operators in the physical domain generalize to the Laplace-Beltrami operator $\Delta$ and the tangential gradient $\nabla$.
The multipatch surface we consider is 
the ellipsoid $\frac{x^2}{3^2} + \frac{y^2}{2^2} + \frac{z^2}{1^2} = 1$ constructed as four maps from unit square reference domains onto four patches that constitute a partition of the ellipsoid. Since this domain has no boundary, we, for a well-posed problem, require that the solution average is zero.
We construct a problem on this domain with a known analytical solution from the ansatz  $u = \sin(4x) \cos(3y)|_{\Omega}$. A numerical solution on this surface is presented in figures~\ref{fig:ellips:a}--\ref{fig:ellips:b}, where we note that the parameterizations are collapsed at two of the poles. The singularities in these maps seem quite nice in comparison to the cusp geometries, so we choose $\gamma=1$ and use $\delta = h^{4\gamma p/(\gamma + 1)}$ as regularization parameter. In figures~\ref{fig:ellips:c}--\ref{fig:ellips:d}, we present convergence results where we note optimal order convergence in both $H^1(\Omega)$-seminorm and in $L^2(\Omega)$-norm using matching approximation spaces of order $p=1,2,3$.

\section{Conclusions}
\label{sec:conclusions}

In this paper, we have presented a generic method for robustly dealing with singular maps in multipatch IGA where the reference domain of each patch is allowed to be trimmed. The main technique is a specific regularization of the Riemannian metric tensor. We want to highlight the following observations:
\begin{itemize}
\item The proposed regularization is very convenient since it can be applied practically without knowledge regarding the behavior of the singular map, and it is applicable to several types of singular parameterizations. For instance, you do not have to know where singular points occur or if the map is actually singular.

\item The cost of the regularization is that you solve the small eigenvalue problem $G a = \lambda a$ at each quadrature point. On the other hand, you no longer have to compute the inverse $G^{-1}$ at each quadrature point. 

\item Also convenient is that the method uses standard isogeometric approximation spaces, where no special precaution is taken close to the problematic parts of the reference patch boundaries, which means that the approximation space may contain functions that do not belong to the correct Sobolev space. Our numerical experiments indicate that the regularized weak form provides sufficient stability for this not to pose a problem.

\item By deriving the method entirely in the reference domains, we naturally get a formulation with correct scalings of the Nitsche penalty term with respect to the Riemannian metric tensors and the mesh sizes in the patches joined at the interface. This is done without defining some specific average operator. 

\item Even if you do not want to regularize, our numerical results show that it is preferable, for robustness and accuracy, to use our method with $\delta=0$ rather than using a naive approach. This is thanks to analytically handling cancellations in the computation of $|G_i|^{1/2} G_i^{-1}$, rather than naively computing this expression as a product of two terms.
\end{itemize}

While the results presented are very promising, the proposed method has so far been thoroughly investigated only for the Poisson problem using our model singular parameterization for the patch geometries. Therefore, we propose the following two future developments:
\begin{itemize}
	\item The regularization procedure was designed with parameterizations in mind where certain dimensions of the parametric geometry diminish as the singularity is approached. Since our numerical experiments only consider the model singular parameterization of this type, a natural progression is to explore other types of singular parameterizations. For instance, the parameterization with colinear map derivatives mentioned at the end of Section~\ref{sec:brief-analysis}. An initial step would be to assess the robustness of the method numerically on a multipatch geometry constructed with such parameterizations.
	\item Given that the regularization is PDE-specific, it would be valuable to generalize it to other problems as outlined in Remark~\ref{rem:other-pde}. Specifically, for surface multipatch geometries, higher-order PDEs such as shell or membrane problems are of interest, where the terms to regularize may include derivatives of the metric.
\end{itemize}

\bibliographystyle{habbrv}
\footnotesize{
\bibliography{ref-singular}
}

\paragraph{Acknowledgements.}
This research was supported in part by the Swedish Research Council Grants Nos.\  2017-03911, 2021-04925,  and the Swedish Research Programme Essence.

\bigskip
\bigskip
\noindent
\footnotesize {\bf Authors' addresses:}

\smallskip
\noindent
Tobias Jonsson,  \quad \hfill \addressumushort\\
{\tt tobias.jonsson@umu.se}

\smallskip
\noindent
Mats G. Larson,  \quad \hfill \addressumushort\\
{\tt mats.larson@umu.se}

\smallskip
\noindent
Karl Larsson,  \quad \hfill \addressumushort\\
{\tt karl.larsson@umu.se}

\end{document}